\documentclass{gen-p-l}
\usepackage{latexsym}
\usepackage{amssymb}
\usepackage{mathrsfs}
\usepackage{rawfonts}
\input{prepictex}
\input{pictex}
\input{postpictex}
\def\Bbb{\mathbb}
\def\be{\begin{equation}}
\def\ee{\end{equation}}
\def\bea{\begin{eqnarray*}}
\def\eea{\end{eqnarray*}}
\newtheorem{main}{Theorem}

\newtheorem{thm}{Theorem}[section]
\newtheorem{lem}[thm]{Lemma}
\newtheorem{prop}[thm]{Proposition}
\newtheorem{cor}[thm]{Corollary}
\newtheorem{defn}[thm]{Definition}

\def\CC{\mathbb C}
\def\RR{\mathbb R}
\def\ZZ{\mathbb Z}
\def\HH{\mathbb H}
\def\PP{\mathbb P}

\def\+{\oplus}

\def\*{^{\ast}}

\def\bcp{{\Bbb C\Bbb P}}

\def\zt{${\mathbb  Z}_2$}

\newenvironment{romenumerate}{\begin{list}{$(\roman{enumi})$}
{\usecounter{enumi}\setlength{\rightmargin}{\leftmargin}}}{\end{list}}
\newcounter{exmp}

\title[Curvature, Optimal Metrics, and $4$-Manifolds]{Curvature Functionals,  
Optimal Metrics, \\ and the Differential Topology of $4$-Manifolds}

\author{Claude LeBrun}
\address{Department of Mathematics\\
SUNY at Stony Brook\\
Stony Brook, NY 11794-3651}
\email{claude@@math.sunysb.edu}
\thanks{Supported 
in part by  NSF grant DMS-0305865.}  
  
\begin{document}

\begin{abstract} This paper investigates the question of which 
smooth compact  $4$-manifolds admit 
 Riemannian metrics
 that  minimize the 
$L^2$-norm of the curvature tensor. Metrics with this property are called 
 {\em optimal}; 
  Einstein metrics and {\em scalar-flat anti-self-dual} 
  metrics provide us with two interesting classes of examples. 
 Using twistor methods, optimal metrics of the second type are constructed 
  on the connected sums $k\overline{\bcp}_2$ for $k > 5$. 
However,  related constructions also show  that 
large classes of  simply connected $4$-manifolds
 {\em do not} admit any  optimal metrics at all. 
Interestingly,  the difference between existence and non-existence
 turns out to delicately depend  on one's choice of smooth structure;  there are 
smooth $4$-manifolds which carry  
  optimal metrics, but  which  are  homeomorphic to infinitely
  many distinct smooth $4$-manifolds on which no optimal metric exists. 
  \end{abstract}

\maketitle
 
 \section{Introduction}
 
 ``Does every smooth compact manifold admit a best metric?''  
 Ren\'e Thom allegedly  
 first posed this 
 na\"{\i}ve but fundamental  question to Marcel Berger \cite{bercent} 
at some point  in the early 1960s. By the early 1990s, it had emerged from the 
 world of informal discussion to find itself in print, as 
  the leading entry on
 one of  S.-T. Yau's 
 celebrated  problem lists \cite{yaulist2}. 
 
 Of course, Thom's question, as formulated above, seems to be  
  less a problem than  a {\em meta-problem}; after all,  
 we are   being asked to find some interpretation of the word ``best'' which will 
 lead to an interesting conclusion. Nevertheless, the question always had some 
  unambiguous content, 
 because Thom had clarified  his question by means of  a paradigmatic example.
 What he of course had in mind was the classical uniformization theorem, which
 tells us  that   every compact $2$-manifold
 carries metrics of constant  curvature.  This paradigm moreover gives us some vital clues
concerning what we ought to look for. First, the definition of ``best metric" should somehow involve
the  Riemannian curvature, and should be invariant under the action of the diffeomorphism
group. Second, one might hope that metrics of constant
sectional curvature, when they exist, would turn out to be the  ``best metrics" on the manifold
in question. 
And third, we should not    expect  our ``best'' metric
to necessarily be absolutely  unique;  a finite-dimensional moduli space of 
``best metrics" would certainly  be  quite acceptable.

If we agree that  the flat metrics are the best metrics on the $n$-dimensional 
torus $T^n=S^1 \times \cdots \times S^1$, then it  seems rather natural to 
 look for   metrics on other manifolds which are ``as flat as possible,'' 
 in the sense that they
minimize some norm of curvature. For example, one might  try to
minimize the $L^p$-norm of the Riemann curvature tensor  for some fixed  $p >1$. However, this is
simply not a sensible problem  for most choices of $p$;  one can typically find a
sequence of metrics for which the $L^p$ norm of curvature tends to 
zero by just multiplying a fixed metric by a suitable sequence of constants. Indeed, there is 
only one value of $p$ for which this trick does {\em not} work: namely, $p=n/2$,
where $n$ is the dimension of the manifold.

Given a smooth compact $n$-manifold $M$,  and  letting  
$${\mathscr G}_M = \{ \mbox{ smooth Riemannian metrics } g \mbox { on }
M \}, $$ we are thus led to  consider 
the functional 
$${\mathcal K}:{\mathscr G}_M {\longrightarrow}  {\RR}  $$
given by 
	$${\mathcal K}(g) =  \int_{M}|{\mathcal R}_{g}|_{g}^{n/2}d\mu_{g}, $$
where ${\mathcal R}$ denotes the Riemann curvature tensor, 
$|{\mathcal R} |$ is its point-wise norm with respect to the metric, and
and 
 $d\mu$ is the  $n$-dimensional  volume measure determined by the metric.
  Berger \cite{bercent} 
 has  suggested the minima of $\mathcal K$
 as  natural  candidates for Thom's  ``best metrics.'' 
Let us codify this proposal:

 \begin{defn}
 Let $M$ be a smooth compact $n$-dimensional manifold, $n\geq 3$.
 A smooth Riemannian metric $g$ on $M$ will be called an 
 {\em optimal metric} if it is an absolute minimizer of the functional $\mathcal K$,
in the sense that 
 $${\mathcal K}(g') \geq {\mathcal K}(g)$$
 for every smooth Riemannian metric $g'$ on $M$. 
 \end{defn}

Notice  that we have defined an optimal 
metric to be a {\em minimum}, not  just a critical point,  of the
functional $\mathcal K$.
This brings into play  a natural diffeomorphism invariant which is
defined even in the absence of  
 an optimal metric:

  \begin{defn}
 For any  smooth compact $n$-dimensional manifold $M$,
 we define ${\mathcal I}_{\mathcal R}(M)$ to be the  
non-negative real number   given by 
 $${\mathcal I}_{\mathcal R}(M) = \inf_{g\in {\mathscr G}_M}{\mathcal K}(g)
 =\inf_g  \int_{M}|{\mathcal R}_{g}|_{g}^{n/2}d\mu_{g}.$$
  \end{defn}

 Thus ${\mathcal I}_{\mathcal R}$ coincides with the number Berger \cite{bercent}
 calls $\min \|{\mathcal R}\|^{n/2}$. 
 Of course,  our definitions have precisely  been chosen  so that  any metric on 
 $g$ on $M$ automatically satisfies
 $${\mathcal K}(g) \geq {\mathcal I}_{\mathcal R}(M),$$
 with equality iff $g$ is an optimal metric.

  While it generally remains  unclear to what extent optimal metrics really 
   represent an appropriate response to 
 Thom's question, the situation in dimension $4$ is rather    
 encouraging. In particular  \cite{bercent,bes}, an Einstein metric 
 on any compact $4$-manifold is optimal. 
However, the converse is by no means true; and the primary purpose of  this
article is  to explore this aspect of the problem. 
Here is what will emerge:

  \begin{itemize}
  \item We will construct infinitely many new examples of non-Einstein optimal 
  metrics on 
   simply connected compact $4$-manifolds. 
  \item We will show that there are many simply connected compact $4$-manifolds 
  which do not admit  optimal metrics. 
  \item We will see that the existence or non-existence of  optimal metrics 
   depends strictly on the diffeotype of a simply connected $4$-manifold;  it  is not 
  determined
  by the homeotype alone.  
  \item We will calculate the invariant ${\mathcal I}_{\mathcal R}$ 
  for many  simply connected $4$-manifolds (some   
  common-garden, others a bit more exotic). 
  \item We will show that the value of ${\mathcal I}_{\mathcal R}$
   depends  strictly on the diffeotype of a $4$-manifold. Different
  differentiable structures on an underlying topological $4$-manifold can
  often be distinguished by the fact that the corresponding values of 
  ${\mathcal I}_{\mathcal R}$ are different. However, 
  we will also see many examples of distinct 
  differentiable structures which {\em cannot} be distinguished in this way. 
    \end{itemize}

\section{Four-Dimensional Geometry}

Our investigation of  optimal metrics on $4$-manifolds will
necessarily presuppose a certain familiarity with  the rudiments of 
$4$-dimensional geometry and topology. The present section 
will attempt to offer a quick introduction to some of this
essential background material.

The notion of self-duality plays a  fundamental r\^ole in  four-dimensional Riemannian geometry. 
If $(M,g)$ is an oriented Riemannian $4$-manifold, the 
Hodge star operator 
$$\star: \Lambda^2 \to \Lambda^2$$
satisfies $\star^2=1$, and so yields a decomposition
\begin{equation} 
\Lambda^2 = \Lambda^+ \oplus \Lambda^- , 
\label{deco} 
\end{equation}
where $\Lambda^+$ is  the $(+1)$-eigenspace of
$\star$, and  $\Lambda^-$ is the $(-1)$-eigenspace. 
Both $\Lambda^+$ and $\Lambda^-$ are  rank-3 vector bundles
over $M$. Reversing the orientation of $M$ interchanges these two bundles. 

\begin{defn}
On any smooth oriented $4$-manifold, 
sections of $\Lambda^+$ are called {\em self-dual 
2-forms}, whereas sections of  $\Lambda^-$ are called {\em anti-self-dual
2-forms}.
\end{defn}

Because  the curvature of any connection is 
a bundle-valued $2$-form, the decomposition  (\ref{deco}) allows one to  break any curvature tensor 
up  into more primitive pieces. This idea has  
 particularly important ramifications when applied to the Riemannian curvature  of the metric 
itself. Indeed, first notice that, by raising an index, the 
Riemann curvature tensor may be reinterpreted as a linear map 
$\Lambda^2 \to \Lambda^2$,  called the {\em curvature operator}. But
decomposing the 2-forms according to  (\ref{deco})
then allows us to view  this linear map  as consisting of four blocks: 
\begin{equation}
\label{curv}
{\mathcal R}=
\left(
\mbox{
\begin{tabular}{c|c}
&\\
$W_++\frac{s}{12}$&$\mathring{r}$\\ &\\
\cline{1-2}&\\
$\mathring{r}$ & $W_-+\frac{s}{12}$\\&\\
\end{tabular}
} \right) . 
\end{equation}
Here $W_\pm$ are the trace-free pieces of the appropriate blocks,
and  are  called the
self-dual and anti-self-dual Weyl curvatures, respectively. 
The scalar curvature  $s$ is understood to act by scalar multiplication,
whereas the   
     trace-free Ricci curvature
$\mathring{r}=r-\frac{s}{4}g$ 
acts on 2-forms by
$$\varphi_{ab} \mapsto ~ 
 \mathring{r}_{ac}{\varphi^c}_{b}-
\mathring{r}_{bc}{\varphi^c}_{a}.$$

An important feature of the decomposition (\ref{deco}) is that it is {\em conformally invariant}, in the sense
that it is unchanged if $g$ is replaced by $u^2g$, where $u$ is an arbitrary smooth positive
function.  Similarly, the self-dual and anti-self-dual Weyl curvatures are 
also conformally invariant (when considered 
as sections of $\Lambda^2 \otimes \mbox{End}(TM)$).

Since our objective is to better understand 
metrics on $4$-manifolds which 
minimize the quadratic curvature integral  
$${\mathcal K}(g)= \int |\mathcal R|^2 d\mu , $$
it is highly  relevant that there are two other quadratic curvature
integrals which actually compute topological invariants. 
Indeed, no matter which  metric $g$ 
we choose on a smooth compact oriented $4$-manifold $M$,  
 the  generalized Gauss-Bonnet theorem \cite{aw} 
tells us  that the Euler characteristic is given by  
\begin{equation}
\label{gb}
\chi (M) = \frac{1}{8\pi^2} \int_M \left(\frac{s^2}{24}+
 |W_+|^2+|W_-|^2-\frac{|\mathring{r}|^2}{2}  \right)d\mu , 
\end{equation}
while  the 
Hirzebruch signature theorem  \cite{hirz} tells us  that 
the {\em signature} is given  by 
\begin{equation}
\label{sig}
\tau (M) = \frac{1}{12\pi^2}\int_M  \left(
 |W_+|^2-|W_-|^2  \right)d\mu ~ .
\end{equation}
Let us recall that  the 
signature of a smooth compact $4$-manifold may be defined in terms of the  
{\em intersection pairing}
\begin{eqnarray*}
\smile  :
H^{2}(M, {\mathbb R})\times H^{2}(M, {\mathbb R})	
 & \longrightarrow & ~~~~ {\mathbb R}  \\
	( ~ [\varphi ] ~ , ~ [\psi ] ~) ~~~~~ & 
	\mapsto  & \int_{M}\varphi \wedge \psi 
\end{eqnarray*}
on de Rham cohomology. By Poincar\'e duality, this
is a non-degenerate pairing; and it is symmetric, since 
$2$-forms commute with respect to the wedge product. 
We may therefore find a basis for $H^2(M, \RR )$ in 
which the intersection pairing is  represented by the  diagonal matrix  	   
	    $$\left[ 
	  \begin{array}{rl}  \underbrace{
	 \begin{array}{ccc}
	   		 1 &  &	  \\
	   		  &	\ddots &   \\
	   		  &	 & 1
	   	  \end{array}}_{b_{+}(M)}
	   	   & 
	   		   \\ 
	  
	 {\scriptstyle b_{-}(M)}  \!
	    \left\{\begin{array}{r}
	    \\
	    \\
	    \\
	    \\
	    \end{array}
	    \right. \! \! \! \! \! \! \! \! \! \! \! \! \! \! \! 
	    &\begin{array}{ccc}
	   		 -1	&  &   \\
	   		  &	\ddots &   \\
	   		  &	 & -1
	   	  \end{array}
	    \end{array} 
	    \right] ,
	    $$ 
 and the numbers $b_\pm (M)$ are then topological invariants of $M$.
Their difference 
$$\tau (M) = b_+ (M) - b_-(M)$$
is the signature of $M$,
whereas their sum 
$$b_2 (M) = b_+(M) + b_-(M)$$
is just the second Betti number. 

A more concrete interpretation of the numbers $b_\pm (M)$ can be given 
by using a bit of Hodge theory. Since   every de Rham class
on $M$ has a unique harmonic representative with respect to $g$, 
we have a canonical identification 
$$H^2(M,{\mathbb R}) =\{ \varphi \in \Gamma (\Lambda^2) ~|~
d\varphi = 0, ~ d\star \varphi =0 \} .$$
But the Hodge star operator $\star$ defines an involution of 
the right-hand side. 
We thus obtain  a direct sum decomposition
\begin{equation}
	H^2(M, {\mathbb R}) = {\mathcal H}^+_{g}\oplus {\mathcal H}^-_{g},
	\label{harm}
\end{equation}
where
$${\mathcal H}^\pm_{g}= \{ \varphi \in \Gamma (\Lambda^\pm) ~|~
d\varphi = 0\} $$
are the spaces of self-dual and anti-self-dual harmonic forms.
The intersection form  is then   positive-definite on  
${\mathcal H}^+_{g}$, and  negative-definite on 
${\mathcal H}^-_{g}$, so we have $$b_{\pm} (M)= \dim {\mathcal H}^\pm_{g}.$$ 
Notice that the spaces ${\mathcal H}^\pm_{g}$ only depend on the
conformal class of the metric. 

One can easily construct $4$-manifolds with any desired values of  $b_+$ and $b_-$
by means of the following construction:

\begin{defn}
Let $M_{1}$ and $M_{2}$ be smooth connected compact oriented $n$-manifolds. 
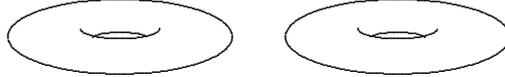
\begin{center}
\mbox{
\beginpicture
\setplotarea x from 0 to 290, y from 0 to 60
\ellipticalarc axes ratio 3:1  360 degrees from 140 40
center at 110 30
\ellipticalarc axes ratio 3:1  -360 degrees from 185 40
center at 215 30
\ellipticalarc axes ratio 4:1 -180 degrees from 125 33
center at 110 33
\ellipticalarc axes ratio 4:1 145 degrees from 120 30
center at 110 29
\ellipticalarc axes ratio 4:1 180 degrees from 200 33
center at 215 33
\ellipticalarc axes ratio 4:1 -145 degrees from 205 30
center at 215 29
\endpicture
}
\end{center}
Their connected sum $M_{1}\# M_{2}$ is then the smooth connected oriented $n$-manifold obtained 
by deleting a small ball from each manifold 
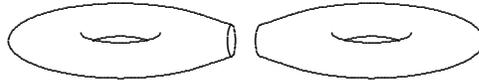
\begin{center}
\mbox{
\beginpicture
\setplotarea x from 0 to 290, y from 0 to 60
\ellipticalarc axes ratio 3:1  270 degrees from 145 40
center at 115 30
\ellipticalarc axes ratio 3:1  -270 degrees from 180 40
center at 210 30
\ellipticalarc axes ratio 4:1 -180 degrees from 130 33
center at 115 33
\ellipticalarc axes ratio 4:1 145 degrees from 125 30
center at 115 29
\ellipticalarc axes ratio 4:1 180 degrees from 195 33
center at 210 33
\ellipticalarc axes ratio 4:1 -145 degrees from 200 30
center at 210 29
\ellipticalarc axes ratio 1:4 360 degrees from 157 36
center at 157 30
\ellipticalarc axes ratio 1:3 180 degrees from 168 36
center at 168 30
{\setlinear 
\plot 145 40        157 36   /
\plot 145 20        157 24  /
\plot 180 40        168 36   /
\plot 180 20        168 24   /
}
\endpicture
}
\end{center}
and identifying the
resulting $S^{n-1}$ boundaries 
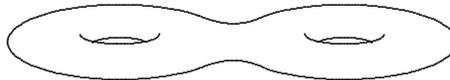
\begin{center}
\mbox{
\beginpicture
\setplotarea x from 0 to 290, y from 0 to 60
\ellipticalarc axes ratio 3:1  270 degrees from 150 40
center at 120 30
\ellipticalarc axes ratio 3:1  -270 degrees from 175 40
center at 205 30
\ellipticalarc axes ratio 4:1 -180 degrees from 135 33
center at 120 33
\ellipticalarc axes ratio 4:1 145 degrees from 130 30
center at 120 29
\ellipticalarc axes ratio 4:1 180 degrees from 190 33
center at 205 33
\ellipticalarc axes ratio 4:1 -145 degrees from 195 30
center at 205 29
{\setquadratic 
\plot 150 40    163 37    175 40   /
\plot 150 20    163 23    175 20  /
}
\endpicture
}
\end{center}
via a reflection. 
\end{defn}

If $M_1$ and $M_2$ are  simply connected $4$-manifolds, then $M= M_1\# M_2$
is also simply connected, and has  $b_\pm (M) = b_\pm (M_1) + b_\pm (M_2)$. 
Now let us use 
$\bcp_{2}$ denote the complex projective plane with its {\em standard} 
orientation, and $\overline{\bcp}_{2}$ denote the same smooth $4$-manifold
with the {\em opposite} orientation. 
Then the 
 iterated connected sum
 $$j \bcp_{2}\# k \overline{\bcp}_{2}
 =
\underbrace{\bcp_{2}\# \cdots \# \bcp_{2}}_{j } 
\# \underbrace{\overline{\bcp}_{2}\# \cdots \# \overline{\bcp}_{2}}_{k }$$
is a simply connected $4$-manifold with $b_+=j$ and $b_-=k$. 
Notice that  $\chi (j \bcp_{2}\# k \overline{\bcp}_{2} )= 
2+j+k$ and  that $\tau (j \bcp_{2}\# k \overline{\bcp}_{2} )= j-k$. 

These $4$-manifolds are {\em non-spin}, meaning that their tangent bundles
have $w_2\neq 0$. For a simply connected compact $4$-manifold $M$, this is equivalent to 
saying that $M$ contains a compact oriented  surface of
odd self-intersection. 

Is this a complete list of the simply connected non-spin $4$-manifolds? Well, yes and no. 
In the affirmative direction,   the remarkable work of 
 Michael Freedman \cite{freedman}  tells us the following:

 \begin{thm}[Freedman]\label{fdmn}
Two smooth compact simply connected oriented $4$-manifolds are orientedly homeomorphic
if and only if 
\begin{itemize}
\item 
they have the same value of $b_+$;
\item 
they have the same value of  $b_-$; and 
\item both are spin, 
or both are 
non-spin.
\end{itemize} 
\end{thm}

Thus, up to {\em homeomorphism}, the connected sums $j \bcp_{2}\# k \overline{\bcp}_{2}$
provide us with a complete list of the simply connected {\em non-spin} $4$-manifolds. 
However, many of these topological $4$-manifolds turn out to have infinitely many
distinct {\em known} smooth structures, and it is generally thought that many of
these manifolds will turn out to have  exotic smooth structures that no one has
even yet imagined.

Freedman's surgical techniques also allow one to classify 
simply connected {\em topological} $4$-manifolds, but the 
classification is much more involved. 
One of the key ingredients that makes it possible to give Theorem \ref{fdmn}
the simple phrasing used above  
is the main theorem of Donaldson's thesis \cite{donaldson}, which showed that 
the anti-self-dual Yang-Mills equations implied previously unsuspected constraints
on the homotopy types of smooth $4$-manifolds:

\begin{thm}[Donaldson]\label{dndsn}
Let $M$ be any smooth  compact simply connected  $4$-manifold with $b_+=0$.
Then $M$ is homotopy equivalent to $S^4$ or to a connected sum $k \overline{\bcp}_2$.
\end{thm}

In particular, if $M$ is a  simply connected differentiable 
$4$-manifold with $b_+(M)=0$ and $b_-(M)\neq 0$, 
 this result tells us that $M$ cannot be spin. 

We have just  seen that 
Theorem \ref{fdmn} allows us to  compile a complete list of simply connected
non-spin homeotypes. In the spin case,  the situation remains more
unsettled, but a conjectural complete list of smoothable simply connected 
spin  homeotypes   consists of $S^4$,
$S^2\times S^2$, $K3$, their connect sums, and orientation reverses
of these. Here  $K3$ means  the 
unique  simply connected smooth compact $4$-manifold admitting 
a complex structure with  $c_1=0$. This $4$-manifold has $b_+=3$ 
and $b_-=19$.  An  interesting model of $K3$  was discovered by 
Kummer, who considered the involution of $T^4$ with $16$ fixed points
which arises as the product of two copies of the 
Weierstrass involution of an elliptic curve:

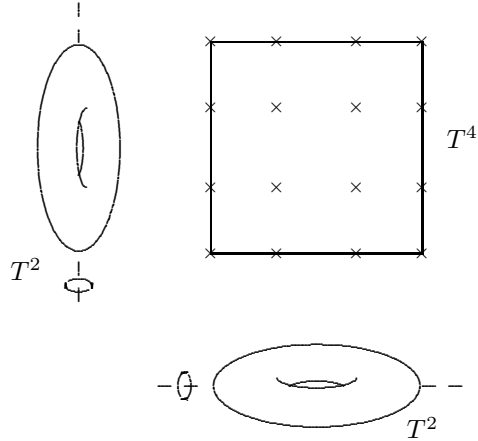
\begin{center}
\mbox{
\beginpicture
\setplotarea x from 0 to 200, y from -5 to 200
\putrectangle corners at 80 160 and 160 80
\put {$_{\times}$} [B1] at 80 105
\put {$_{\times}$} [B1] at 105 105
\put {$_{\times}$} [B1] at 135 105
\put {$_{\times}$} [B1] at 160 105

\put {$_{\times}$} [B1] at 80 80
\put {$_{\times}$} [B1] at 105 80
\put {$_{\times}$} [B1] at 135 80
\put {$_{\times}$} [B1] at 160 80

\put {$_{\times}$} [B1] at 80 135
\put {$_{\times}$} [B1] at 105 135
\put {$_{\times}$} [B1] at 135 135
\put {$_{\times}$} [B1] at 160 135

\put {$_{\times}$} [B1] at 80 160
\put {$_{\times}$} [B1] at 105 160
\put {$_{\times}$} [B1] at 135 160
\put {$_{\times}$} [B1] at 160 160

\put {$T^2$} [B1] at 160 10 
\put {$T^2$} [B1] at  10 70
\put {$T^4$} [B1] at  175  120
\ellipticalarc axes ratio 5:2  360 degrees from 150 40
center at 120 30
\ellipticalarc axes ratio 4:1 -180 degrees from 135 33
center at 120 33
\ellipticalarc axes ratio 4:1 145 degrees from 130 30
center at 120 29
\ellipticalarc axes ratio 2:5  -360 degrees from 40 150 
center at 30 120 
\ellipticalarc axes ratio 1:4  180 degrees from  33 135
center at  33 120
\ellipticalarc axes ratio 1:4 -145 degrees from  30 130
center at  29 120
\ellipticalarc axes ratio 1:2 140 degrees from 70 35
center at 70 30
\ellipticalarc axes ratio 1:2 140 degrees from 70 25
center at 70 30
\ellipticalarc axes ratio 2:1 140 degrees from 35 68
center at 30 68
\ellipticalarc axes ratio 2:1 140 degrees from 25 68
center at 30 68
\arrow <2pt> [.1,.3] from 70 35 to 71 35
\arrow <2pt> [.1,.3] from 70 25 to 69 25
\arrow <2pt> [.1,.3] from 35 68 to 35 67
\arrow <2pt> [.1,.3] from 25 68 to 25 69
{\setlinear
\setdashes 
\plot   30 62  30 82  /
\plot  60 30 80 30   /
\plot  160 30 180 30   /
\plot   30 160  30  180  /
}
\endpicture
}
\end{center}
Kummer's model of $K3$ is then obtained from the orbifold
$T^4/\ZZ_2$ by replacing each singular point with a
$\bcp_1$ of self-intersection $-2$. Analogous
 constructions  will turn out to play a central  r\^ole 
in this paper. 

\section{Optimal Geometries in Dimension Four}

A Riemannian metric is said to be {\em Einstein}
 if it has constant Ricci curvature.
 Since the Ricci curvature of $g$ is by definition the function 
  $v\mapsto r (v,v)$
  on the unit tangent bundle $g(v,v)=1$, where $r$ denotes the Ricci tensor, 
   this is clearly equivalent to the requirement  that 
 $$r = \lambda g$$
 for some  constant $\lambda$.
This in turn can be rewritten as the pair
 of conditions
 $$
\mathring{r}  =  0 , ~~~~~
s  =  \mbox{constant} ,
$$
where $s= \mbox{trace}_g(r)$ is the scalar curvature, and where $\mathring{r}=r-\frac{s}{n}g$ is the 
trace-free part of the Ricci tensor. However, double contraction of the second Bianchi
identity tells  us that 
 $$\nabla\cdot \mathring{r}= {\textstyle \frac{n-2}{n}}~ds,$$
 so, in any dimension $n\neq 2$, a metric is Einstein iff it satisfies the equation 
 $$\mathring{r}=0.$$
 
 An important motivation for the study of optimal metrics is that \cite{bercent,bes} any Einstein
 metric is optimal in dimension $4$. Indeed, let $M$ 
 be a smooth compact $4$-manifold, and let $g$ be an arbitrary Riemannian metric
 on $M$. Then the  $4$-dimensional Gauss-Bonnet formula (\ref{gb})  allows us
 to rewrite
 $$
 {\mathcal K}(g) = \int_M |{\mathcal R}|^2d\mu_g
 $$
 as 
 \begin{equation}
\label{chess}
{\mathcal K}(g) = 8\pi^2 \chi (M) + \int_M |\mathring{r}_g|^2d\mu_g ~. 
\end{equation}
Thus any metric $g$ with  $\mathring{r}=0$ minimizes $\mathcal K$;
and when such a metric exists, the Einstein metrics are the {\em only}
optimal metrics on $M$. 
 
 However, similar arguments also show that smooth compact $4$-manifolds 
 often do not admit Einstein metrics, even in the simply connected case.
 Indeed, a judicious combination of (\ref{gb}) with (\ref{sig}) reads
 \begin{equation}
\label{pnumb}
(2\chi + 3\tau )(M) = \frac{1}{4\pi^2}\int_M\left( \frac{s^2}{24}+2|W_+|^2-
 \frac{|\mathring{r}|^2}{2}\right)d\mu_g,
\end{equation}
 so that a compact oriented $4$-manifold $M$ can only admit an Einstein metric if it satisfies the
 {\em Hitchin-Thorpe inequality} \cite{hit,tho} 
 $$(2\chi + 3\tau )(M)\geq 0.$$
Thus, for example, the simply connected $4$-manifolds
$k\overline{\bcp}_2$ do not admit an Einstein metric when $k> 4$,
since these spaces have $2\chi+3\tau = 4-k$.

 In the boundary case of the Hitchin-Thorpe inequality, we also get
 a striking amount of additional information. Indeed, 
 if $(2\chi + 3\tau )(M)=0$   and $M$ admits an Einstein metric $g$, then 
 the  Riemannian connection on $\Lambda^+\to M$ 
 must be flat, since its curvature tensor is algebraically
determined by $r$ and $W_+$; such a metric is said to be 
{\em locally hyper-K\"ahler}. Any locally hyper-K\"ahler
$4$-manifold is finitely covered 
 \cite{bes,hit}
by a flat $4$-torus
or a Calabi-Yau  $K3$. In particular, 
a simply connected $4$-manifold with  $2\chi + 3\tau =0$
can admit an Einstein metric only if it is diffeomorphic to $K3$. 
In particular, $4\overline{\bcp}_2$ does not admit an Einstein
metric, since it has $2\chi+3\tau =0$, but is not even homotopy
equivalent to $K3$.

 It should therefore seem a bit reassuring that there are many 
 simply connected $4$-manifolds which carry optimal metrics which 
 are not Einstein \cite{lebicm}. To see this, we may begin by  observing   
 that (\ref{gb}) and (\ref{sig}) allow one to re-express ${\mathcal K}$
 as \begin{equation}
\label{knut}
{\mathcal K}(g) = -8\pi^2 (\chi + 3\tau )(M) + 2\int_M \left(\frac{s^2}{24}+2|W_+|^2\right) d\mu_g ~. 
\end{equation}
This brings another important class of  metrics to the fore:

\begin{defn}
If $M$ is a smooth oriented $4$-manifold, a Riemannian metric $g$
on $M$ will be said to be {\em anti-self-dual} (or, for brevity,  {\em ASD}) 
if its self-dual Weyl curvature is identically zero: 
$$W_+\equiv 0.$$
A metric $g$ will be called 
{\em scalar-flat}  (or, more briefly, {\em SF}) 
if it satisfies
$$s=0.$$
Finally, we will say that $g$ is {\em scalar-flat anti-self-dual} (or {\em SFASD}) if it satisfies both of 
these conditions. 
\end{defn}

Equation (\ref{knut}) then immediately yields the following:
\begin{prop}
Suppose that $M$ is a smooth compact oriented $4$-manifold. If $M$ carries
a scalar-flat anti-self-dual metric $g$, then $g$ is optimal. When this happens,  moreover, 
every other optimal metric $g^\prime$ on $M$ is SFASD, too. 
\end{prop}

This fact seems to have first been noticed by Lafontaine \cite{laf}, who simultaneously   
discovered an important
 obstruction to the existence of  SFASD metrics.
 Indeed, equation  (\ref{pnumb})  and our previous 
  discussion of the locally hyper-K\"ahler manifolds 
immediately gives us the following 
  upside-down version of the Hitchin-Thorpe 
inequality: 

\begin{prop}[Lafontaine]\label{fine}
Let $(M,g)$ be a  compact scalar-flat anti-self-dual $4$-manifold. 
Then 
$$(2\chi + 3 \tau ) (M) \leq 0,$$
with equality if and only if $(M,g)$ is finitely covered by a flat $4$-torus or  a Calabi-Yau $K3$. 
\end{prop}

However, a completely different set of topological 
constraints is imposed by the  following result \cite{lsd}:

\begin{prop}\label{form}
Let $(M,g)$ be a   compact scalar-flat anti-self-dual $4$-manifold. 
Then either 
\begin{itemize}
\item $b_+(M)=0$; or 
\item $b_+(M)=1$, and $g$ is a scalar-flat K\"ahler metric; or else 
\item $b_+(M)=3$, and $g$ is a  hyper-K\"ahler metric. 
\end{itemize}
\end{prop}
 
 \begin{proof} Recall that $b_+(M)$ is exactly the dimension of the 
 space of harmonic self-dual $2$-forms on $M$. However, any self-dual $2$-form
 $\varphi$ on any Riemannian $4$-manifold satisfies the Weitzenb\"ock formula
 \cite{bourg} 
 $$(d+d^*)^2 \varphi = \nabla^*\nabla \varphi - 2 W_+(\varphi , \cdot ) + \frac{s}{3}\varphi , $$
 so if $\varphi$ is harmonic and if $g$ is SFASD we obtain 
 $$0= \int_M \langle \varphi , \nabla^*\nabla \varphi \rangle d\mu = \int_M |\nabla \varphi |^2 d\mu.$$
 Thus, when $(M,g)$ is a compact SFASD manifold, 
 $b_+(M)$ is exactly  the dimension of the space of  the parallel self-dual $2$-forms.
 
 Now $SO(4)$ is a double cover of $SO(3)\times SO(3)$, where the factor projections to   
 $SO(3)$ are given by  its action on $\Lambda^\pm$. The subgroup of $SO(4)$ stabilizing a
 non-zero element of $\Lambda^+$ is thus the double cover $U(2)$ of 
 $SO(2)\times SO(3)$, whereas the subgroup acting trivially on $\Lambda^+$, or even
 on a $2$-dimensional subspace of it,  is 
 the universal cover $SU(2)=Sp(1)$  of $SO(3)$. Thus an oriented Riemannian $4$-manifold 
 with a non-trivial parallel self-dual $2$-form has holonomy $\subset U(2)$, and is K\"ahler,
 whereas the existence of $2$  independent parallel self-dual $2$-forms would force
 the manifold to have holonomy $\subset Sp(1)$, and so to be hyper-K\"ahler. 
   \end{proof}

   Let us now assume that $(M,g)$ is a {\em simply connected} SFASD manifold, 
   and see what the above results now tell us. By Proposition \ref{form}, the 
   only possibilities for $b_+(M)$ are $0,1$ and $3$. If $b_+=3$,
   we would have a simply connected    hyper-K\"ahler manifold; 
   such an object is necessarily a $K3$ surface. If $b_+=1$, we have
   a simply connected complex surface with non-trivial canonical line bundle and
   a K\"ahler metric of zero scalar curvature;
   by a plurigenus vanishing theorem of Yau \cite{yauruled} 
   and the Enriques-Kodaira classification \cite{bpv}, 
   such a complex surface must be obtained from $\bcp_2$ by blowing up and down, and hence diffeomorphic to either 
   $S^2\times S^2$ or a connected sum $\bcp_2 \# k \overline{\bcp}_2$; and since 
  Proposition \ref{fine} tells us  that $2\chi + 3\tau = 4 -4b_1+5b_+ -b_-= 9-b_-$ is negative,
  we conclude in this case that $M$ is diffeomorphic to $\bcp_2 \# k \overline{\bcp}_2$ 
  for some $k> 9$. Finally, if $b_+=0$, 
  Theorems \ref{dndsn} and \ref{fdmn} 
  tell us that $M$ is at least {\em homeomorphic} to either $S^4$ or a connected sum 
  $k\overline{\bcp}_2$; and since 
    Proposition \ref{fine} tells us  that $2\chi + 3\tau = 4 -b_-$ is negative, we conclude
    in this case that $M$ is homeomorphic to $k\overline{\bcp}_2$ for some $k>4$.
    Summarizing, we have \cite{lsd}

\begin{prop}\label{class} 
Let $M$ be a smooth compact 
 {\em simply connected}  $4$-manifold. 
 If $M$ admits a scalar-flat anti-self-dual metric $g$, then 
\begin{itemize}
\item $M$ is homeomorphic to $k \overline{\bcp}_2$ for some $k \geq 5$; or 
\item $M$ is diffeomorphic to $\bcp_2\#k \overline{\bcp}_2$ for some  $k \geq 10$; or else 
\item $M$ is diffeomorphic to $K3$.
\end{itemize}
\end{prop}

\bigskip 
 
 \noindent 
 A major objective of this paper is to prove  the following partial 
 converse: 
 
 \bigskip 
 
 \begin{main}\label{able}
 A   simply connected $4$-manifold $M$  admits scalar-flat anti-self-dual 
 metrics if 
 \begin{itemize}
 \item  $M$ is diffeomorphic to  $k \overline{\bcp}_2$ for  some $k  \geq 6$; or 
\item $M$ is diffeomorphic to $\bcp_2\#k \overline{\bcp}_2$ for some  $k \geq 14$; or  
\item $M$ is diffeomorphic to $K3$.
\end{itemize}
In particular, each of these manifolds carries optimal metrics. 
 \end{main}
 
 It is worth emphasizing  that,  except in the $K3$ case,  the   
 optimal metrics of Theorem \ref{able} are necessarily  non-Einstein. 
 
 On the other hand, Corollary \ref{class} and a computation of ${\mathcal I}_{\mathcal R}$
 will allow us to show 
  that the existence of  optimal metrics 
 is highly sensitive to the choice of smooth structure: 
 
 \begin{main}\label{baker}
For each $k \geq 9$, the
topological manifold $\bcp_2\#k \overline{\bcp}_2$ admits 
 infinitely many ``exotic'' smooth structures for which the corresponding
 smooth compact $4$-manifold does not admit optimal metrics. Similarly, 
  the
topological manifold $K3$ 
admits 
 infinitely many exotic smooth structures for which the corresponding
 smooth $4$-manifold does not admit optimal metrics.  
 \end{main}

  Similar ideas will also allow us to prove the non-existence of 
  optimal metrics for smooth manifolds representing many more 
   homeotypes:

 \begin{main}\label{charlie}
 If $j\geq 2$ and   $k\geq 9j$, 
 the smooth simply connected $4$-manifold
 $j \bcp_2\#k \overline{\bcp}_2$ 
 does not admit optimal metrics.
 Moreover,  if $j\geq 5$ and $j\not\equiv 0\bmod 8$, the
 underlying topological manifold of this space admits infinitely
 many distinct differentiable structures for which no
 optimal metric exists. 
 \end{main}

\section{Constructing anti-self-dual metrics}
\label{construct}

 The condition of anti-self-duality is {\em  conformally invariant}; if $g$ is an ASD
 metric, so is $u^2g$, for any $u > 0$. 
The strategy of our proof of Theorem \ref{able} will be to 
first construct a family of anti-self-dual conformal classes of metrics on $k\overline{\bcp}_2$, 
$k\geq 6$,
and then show that some of the constructed conformal classes contain 
scalar-flat metrics. Our approach to both aspects of this problem will be carried out
using methods of complex analysis via the Penrose 
 {\em twistor correspondence}, to which we now provide a brief introduction. 

 Given any  
 oriented Riemannian 4-manifold  $(M,g)$, one can construct
an associated
almost-complex $6$-manifold $(Z, J)$,
where   $\pi : Z\to M$  is the 
$S^2$-bundle $S(\Lambda^+)$ of unit self-dual 
2-forms. The almost-complex structure  $J:TZ\to TZ$ preserves the decomposition 
of $TZ$ into horizontal and vertical components with 
respect to the Levi-Civita connection.
On the tangent spaces of each fiber $S^2$, $J$ simply acts 
by rotation by $-90^{\circ}$.  Meanwhile, in the horizontal sub-bundle,
 which we identify 
 with $\pi^*TM$,
$J$ acts 
at $\phi\in S(\Lambda^+)$
 by  $v\mapsto \sqrt{2}(v \lrcorner \phi )^\sharp$.  
 Each fiber $S^2$ of $S(\Lambda^+)\to M$ is thus a $J$-holomorphic curve, and 
 the fiber-wise antipodal map $\sigma: S(\Lambda^+)\to S(\Lambda^+)$
 is $J$-anti-holomorphic, in the sense that $\sigma_*\circ J= -J\circ \sigma_*$. 
A remarkable and non-obvious feature of this construction is that   
the almost-complex structure $J$ is   
actually {\em conformally invariant},  despite the fact that  replacing $g$ with $u^2g$ 
alters the horizontal subspaces on 
$Z=(\Lambda^+-0)/{\mathbb  R}^+$.

 Now recall that an almost-complex manifold is 
 a complex manifold iff it admits sufficiently many  local  
 holomorphic functions.  
   In general, the obstruction \cite{newnir}
to the existence of such  functions is 
  the 
{\em Nijenhuis tensor},  but in the present  case
the Nijenhuis tensor of $(Z,J)$ just amounts to  the self-dual Weyl 
   curvature  $W_+$ of 
 $(M,g)$. 
When $(M,g)$ is anti-self-dual,  $(Z,J)$ thus acquires the structure of a complex manifold  \cite{AHS,pnlg}: 

\begin{thm}[Penrose/Atiyah-Hitchin-Singer] The almost-complex
manifold $(Z,J)$ is a complex 3-manifold iff $W_+=0$. 
Moreover, a complex 3-manifold arises by this construction
iff it admits a fixed-point-free anti-holomorphic 
involution $\sigma: Z\to Z$ and a foliation by $\sigma$-invariant
rational curves ${\bcp}_1$, each of which has normal 
bundle ${\mathcal  O}(1)\oplus {\mathcal  O}(1)$. Finally,
the complex manifold $(Z,J)$ and the real structure 
$\sigma$  suffice  to determine the metric
$g$ on $M$ up to conformal rescaling. \label{pahs}
\end{thm}

\begin{defn} The complex 3-manifold $(Z,J)$
associated with 
 an  anti-self-dual 4-manifold $(M,g)$ by
Theorem \ref{pahs} is called the {\em
twistor space} of $(M,[g])$. 
\end{defn}
   
   \begin{defn} Let $(M,g)$ be an anti-self-dual $4$-manifold,
   let $(Z,J)$ be its twistor space, 
   and let $x\in M$. Then the  holomorphic curve 
    $P_x\subset Z$  
   given by $\pi^{-1}(x)= S(\Lambda^+_x)$ 
     will be called the {\em real twistor line}
    corresponding to $x$. 
   \end{defn}
   
   The moduli space of holomorphic  curves $\bcp_1\subset Z$ near the real 
   twistor lines   is a complex $4$-manifold ${\mathcal M}$, and  is a 
   complexification of the original real $4$-manifold $M$. The term {\em complex twistor
   line} (or just {\em twistor line}) is used to refer to any $\bcp_1$ in this larger family.

One of the cornerstones of the theory of anti-self-dual manifolds is the connected sum
construction of Donaldson and Friedman \cite{DF}. If $M_1$ and $M_2$ admit
anti-self-dual metrics, this allows one to  construct anti-self-dual metrics on the connected sum 
$M_1\# M_2$, provided the  twistor spaces $Z_1$ and $Z_2$ of the given manifolds
satisfy $H^2(Z_j , {\mathcal O}(TZ_j))=0$, $j=1,2$.  An orbifold  generalization of this
construction was later developed by the present author in collaboration with Michael
Singer \cite{lebsing2}, and allows one to build up non-singular anti-self-dual manifolds
by gluing special orbifolds across $\RR\PP^3 \times \RR$ necks. We will now
review those features of this generalized  construction which will be needed in
what follows. 

Let $(Y,g)$ be a compact anti-self-dual manifold with twistor space $Z$, and
let us assume from the outset that $H^2(Z, {\mathcal O}(TZ))=0$. 
Let $\phi : Y \to Y$ be an isometry of $(Y,g)$ with $\phi^2=\text{id}_Y$, and 
assume that $\phi$ has exactly $k$ fixed points, for some 
positive integer $k$. Each fixed point is therefore isolated, and  
in geodesic normal coordinates around any fixed point,  $\phi$ therefore just 
becomes the involution $\vec{v}\to -\vec{v}$ of $\RR^4$; in particular, 
$\phi$ is   orientation-preserving. 
We thus have an induced map ${\phi}^*: S(\Lambda^+)\to S(\Lambda^+)$, 
and this map may be viewed as a holomorphic involution $\hat{\phi} : Z\to Z$. 
The fixed point set of $\hat{\phi}$ then consists of a disjoint union of $k$ real 
twistor lines, one for each fixed point of $\phi$. Let $\tilde{Z}$ be obtained 
by blowing up $Z$ along these $k$ twistor lines, and notice that 
$\hat{\phi}$ induces a holomorphic involution $\tilde{\phi}$ of 
$\tilde{Z}$. The fixed point set of $\tilde{\phi}$ then consists of $k$ quadrics
$\bcp_1\times \bcp_1$ with normal bundle ${\mathcal O}(1,-1)$, and the 
$\tilde{\phi}$ acts on their normal bundles by multiplication by $-1$.
The quotient $\tilde{Z}/\ZZ_2$ is therefore a non-singular compact complex
$3$-fold containing $k$ hypersurfaces $Q_1, \ldots Q_k$, each biholomorphic to 
$\bcp_1\times \bcp_1$, and each with normal bundle ${\mathcal O}(2,-2)$.
The complement of these hypersurfaces is just the twistor space of 
$[Y-\{ \text{fixed points}\}]/\ZZ_2$. In this complement, choose $\ell \geq 0$ 
twistor lines, and blow them up to obtain  $\ell$ hypersurfaces 
$Q_{k+1}, \ldots , Q_{k+\ell }$, each biholomorphic to $\bcp_1\times \bcp_1$,
and each with normal bundle  ${\mathcal O}(1,-1)$.  Let $Z_+$ denote 
this blow-up of $\tilde{Z}/\ZZ_2$, and notice that our original anti-holomorphic
involution $\sigma$ of $Z$ induces an anti-holomorphic involution 
$\sigma_+: Z_+\to Z_+$. 
 
 The next ingredient we will need is a compactification of the twistor space of 
 the Eguchi-Hanson metric. The usual Eguchi-Hanson metric \cite{EH} is a locally asymptotically
 flat hyper-K\"ahler metric on $T^*S^2$ which, up to homothety, is
    the metric-space  completion of the 
 Riemannian metric 
 $$g_{EH, \epsilon}= \frac{d\varrho^2}{1-{\varrho}^{-4}}+ \varrho^2\left( \sigma_1^2+\sigma_2^2+ 
\left[1-{\varrho}^{-4}\right]\sigma_3^2\right)$$
 on $(1, \infty ) \times S^3/\ZZ_2$, where $\{ \sigma_j\}$ is 
the standard  left-invariant co-frame on $S^3/\ZZ_2 = SO(3)$.
However, because this metric is asymptotic to the flat metric on $(\RR^4-\{ 0\})/\ZZ_2$ 
as $\varrho \to \infty$, 
  its conformal class naturally extends
 to an orbifold  ASD conformal metric on $T^*S^2\cup \{\infty\}$, where the added
 point is  singular,  with a neighborhood modeled on $\RR^4/\ZZ_2$. Blowing up the 
 twistor line of this added ``point at infinity''  then yields a non-singular complex $3$-fold
 $\tilde{Z}_{EH}$ which contains a hypersurface $Q\cong \bcp_1\times \bcp_1$ with normal
 bundle ${\mathcal O}(2,-2)$ arising as the exceptional divisor of the blow-up.
The complex $3$-fold $\tilde{Z}_{EH}$ is \cite{lebsing2} a 
small resolution of  the hypersurface 
$$
xy=z^2-t^2\zeta_1^2\zeta_2^2
$$
 in the $\bcp_3$-bundle  
 ${\PP}({\mathcal O}(2)^{\oplus 3}\oplus {\mathcal O})$
 over $\bcp_1$, where $x,y,z\in {\mathcal O}(2)$, $t\in {\mathcal O}$, and $[\zeta_1: \zeta_2]$
 are the homogeneous coordinates on  $\bcp_1$. The small resolutions  replace
  the two singular points $x=y=z=0$ with rational curves.
 
Finally, consider the Fubini-Study metric on $\bcp_2$, which, up to homothety,  may be 
characterized as  the
the unique $SU(3)$-invariant metric on the complex projective plane. 
Because the isotropy subgroup $U(2)\subset SU(3)$ is so large, 
the Fubini-Study metric has  $W_-=0$ for representation-theoretic reasons;
reversing the orientation, the Fubini-Study metric thus becomes an ASD metric 
on $\overline{\bcp}_2$. Let $Z_{FS}$ denote the twistor space of this metric, 
and let $\tilde{Z}_{FS}$ denote its blow-up along a real twistor line. Explicitly,
$Z_{FS}$ may be realized \cite{AHS} as 
$$\{ ([x_1:x_2:x_3], [y_1:y_2:y_3])\in \bcp_2 \times \bcp_2 ~|~ \sum_jx_jy_j =0\},$$
and we may take the relevant twistor line to be given by $x_1=y_1=0$. 
Blowing up of this twistor line provides us with a preferred hypersurface
$\bcp_1\times \bcp_1$ with normal bundle ${\mathcal O}(1,-1)$. 

Now let $Z_-$ be the disjoint union of $k$ copies of  $\tilde{Z}_{EH}$
and $\ell$ copies of $\tilde{Z}_{FS}$, and let $\sigma_-$ be the real structure
it inherits from the twistor spaces of the Eguchi-Hanson and Fubini-Study metrics. 
 Let $Q_-\subset Z_-$ be the disjoint union 
of $k+\ell$ copies of $\bcp_1\times \bcp_1$, each being a copy of the constructed exceptional
 divisor in a copy of  $\tilde{Z}_{EH}$ or $\tilde{Z}_{FS}$. 
Remembering that $Z_+$ also contains a disjoint union 
 $Q_+=Q_1\sqcup \cdots \sqcup Q_{k+\ell}$ of 
 the same number of copies of $\bcp_1\times \bcp_1$,  we may thus form
 a $3$-dimensional 
complex-analytic space $Z_0$ with normal crossing singularities by 
identifying $Q_+$ with $Q_-$. However, we carry out this identification
according to a few simple rules.
 First of all,  the quadrics $Q_1, \ldots, Q_k$ are each to be identified 
with a quadric in a copy of $\tilde{Z}_{EH}$, while the remaining
quadrics $Q_{k+1} , \ldots Q_{k+\ell}$ are each to be  identified with a quadric
in a copy of $\tilde{Z}_{FS}$. 
 Secondly, 
we always  interchange the factors of $\bcp_1\times\bcp_1$ when gluing 
$Q_+$ to $Q_-$, thereby making the normal bundles of each 
quadric relative to $Z_+$ an $Z_-$  dual to each other.
Finally, we always make our identifications
in such a way that the real structures  $\sigma_+$ and $\sigma_-$ 
agree on the locus $Q$ obtained by identifying $Q_+$ with 
$Q_-$, so that $Z_0$ comes equipped with  an anti-holomorphic 
involution $\sigma_0 : Z_0 \to Z_0$. 

Because we have assumed that $H^2(Z, {\mathcal O}(TZ))=0$, 
one can show \cite{lebsing2} that  $\mbox{\bf Ext}^2_{Z_0}(\Omega^1, {\mathcal O})=0$,
and a generalization of Kodaira-Spencer theory \cite{friedman} then 
yields a versal deformation of $Z_0$,
parameterized by a neighborhood of the the origin in 
$\mbox{\bf Ext}^1_{Z_0}(\Omega^1, {\mathcal O})$.
The generic fiber of this family is non-singular
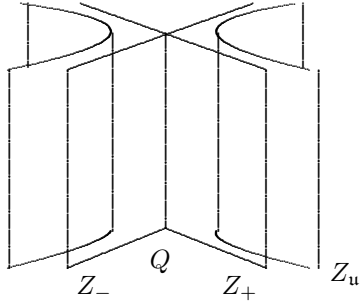
\begin{center}
\mbox{
\beginpicture
\setplotarea x from 0 to 200, y from -20 to 110 
\put {$Z_+$} [B1] at 125  -10
\put {$Z_-$} [B1] at  70  -10
\put {$Z_{\mathfrak u}$} [B1] at 165  -5
\put {$Q$} [B1] at    95 0 
{\setlinear 
\plot 65 100 135 75 /
\plot 60 75  130 100 /
\plot 135 75 135 0 /
\plot 60 75 60 0 /
\plot 60 0 97 15 /
\plot 97 15 135 0 /
\plot 97 89 97 15 /
\plot 117 89 117  14 /
\plot 155 75 155 0 /
\plot 149 100 149 76 /
\plot 38 75 38 0 /
\plot 77 89 77 14 /
\plot 45 100 45 76 /
}
{\setquadratic
\plot  148  75   113  89   143  100 /
\plot  148  0    120 8        113  14 /
\plot 34 75 73 89 39 100 /
\plot 34 0     62 8         73 14  /
}
\endpicture
}
\end{center}
and the real structure $\sigma_0$ extends  to act on the total space
of this family. 

Rather than working with the entire versal family, it is convenient to 
restrict ones attention to certain subfamilies, called {\em standard deformations}.
A {\em $1$-parameter standard deformation} of $Z_0$ is by definition a flat proper
holomorphic map $\varpi : {\mathcal Z}\to {\mathcal U}$ together with
 an anti-holomorphic involution
$\sigma : {\mathcal Z}\to {\mathcal Z}$, such that 
\begin{itemize}
\item
${\mathcal Z}$ is a complex $4$-manifold;
\item  ${\mathcal U}\subset \CC$ is
an open neighborhood of $0$;
\item  $\varpi^{-1}(0)=Z_0$;
\item $\sigma |_{Z_0}= \sigma_0$; 
\item $\sigma$ covers complex conjugation in $\mathcal U$;
\item $\varpi$ is a submersion away from $Q\subset Z_0$; and 
\item near any point of $Q$, $\varpi$ is modeled on $(x,y,z,w)\mapsto xy$. 
\end{itemize}  
When ${\mathfrak u}\in {\mathcal U}\subset \CC$ is  real, non-zero,
and sufficiently small, the corresponding fiber $Z_{\mathfrak u}= \varpi^{-1}({\mathfrak u})$
is a twistor space, and we obtain the following result  \cite{DF,lebsing2}:

\begin{thm} \label{thmb} 
Let $(Y,g)$ be a compact anti-self-dual $4$-manifold equipped with 
an isometric \zt-action with exactly $\ell$ fixed points, for some positive integer $\ell$. 
Let $Z$ denote the twistor space of $(Y, g)$,  and suppose that 
$H^2(Z, {\mathcal O}(TZ))=0$. 
Let $X= Y/\ZZ_2$, and let
 $\tilde{X}$ be the oriented manifold obtained by replacing each 
 singularity of $X$ with a $2$-sphere of self-intersection $-2$. 
 Then, for any integer $k\geq 0$, there are anti-self-dual conformal classes on 
 $\tilde{X}\# k\overline{\bcp}_2$ whose twistor spaces arise as fibers in a 
 $1$-parameter standard deformation of $Z_0=(\tilde{Z}/\ZZ_2)\cup \ell \tilde{Z}_{EH}\cup
 k \tilde{Z}_{FS}$. 
\end{thm}

Similarly, one can define standard deformations depending on several parameters.
For example, a  {\em $2$-parameter standard deformation} of $Z_0$ is by definition a flat proper
holomorphic map $\varpi : {\mathcal Z}\to {\mathcal U}$ together with
 an anti-holomorphic involution
$\sigma : {\mathcal Z}\to {\mathcal Z}$, such that 
\begin{itemize}
\item
${\mathcal Z}$ is a complex $5$-manifold;
\item  ${\mathcal U}\subset \CC^2$ is
an open neighborhood of $(0,0)$;
\item  $\varpi^{-1}(0,0)=Z_0$;
\item $\sigma |_{Z_0}= \sigma_0$; 
\item $\sigma$ covers complex conjugation in $\mathcal U$;
\item $\varpi$ is a submersion away from $Q\subset Z_0$; and 
\item near any point of $Q$, $\varpi$ is modeled on $(x,y,z,v,w)\mapsto (xy,z)$. 
\end{itemize}  
When $({\mathfrak u}_1,{\mathfrak u}_2)\in {\mathcal U}\subset \CC^2$ is  real, 
and sufficiently close to $(0,0)$, with 
${\mathfrak u}_1\neq 0$, the corresponding fiber $Z_{\mathfrak u}= \varpi^{-1}({\mathfrak u}_1,
{\mathfrak u}_2)$
is a twistor space. By extracting such standard deformations from the
 versal deformation of $Z_0$, the same proof tells one the following:

\begin{thm} \label{thumbo} 
Let $Y$ be a compact real-analytic oriented $4$-manifold equipped with 
a \zt-action with exactly $\ell$ fixed points, for some positive integer $\ell$. 
Let $[g_t]$, $t\in (-\varepsilon , 
\varepsilon )$, 
be a 
real-analytic $1$-parameter family of \zt-invariant anti-self-dual conformal metrics
on $Y$.  
 Let $X= Y/\ZZ_2$, and let
 $\tilde{X}$ be the oriented manifold obtained by replacing each 
 singularity of $X$ with a $2$-sphere of self-intersection $-2$. 
Let $Z$ denote the twistor space of $(Y, [g_0])$, 
and suppose that $H^2(Z, {\mathcal O}(TZ))=0$.  
 Then, for any  $k\geq 0$, there is a $2$-parameter standard deformation of 
 $Z_0= (\tilde{Z}/\ZZ_2)\cup \ell \tilde{Z}_{EH}\cup
 k \tilde{Z}_{FS}$ such that, for all small real real numbers   ${\mathfrak u}_1 > 0$
 and ${\mathfrak u}_2$, the fiber $Z_{({\mathfrak u}_1, {\mathfrak u}_2)}$
 is the twistor spaces of a family of an anti-self-dual metric on 
  $\tilde{X}\# k\overline{\bcp}_2$, and such that, for all sufficiently small real numbers
  $t$,  the fiber over $Z_{(0,t)}$ is the complex-analytic space with normal crossings built from 
  the twistor space of $(Y,[g_t])$ by analogy to the construction of 
  $Z_0$. 
\end{thm}

\section{Conformal Green's Functions}
\label{greatsmall}

Let $(M,g)$ be a compact Riemannian $4$-manifold,
 and assume that 
its Yamabe Laplacian $\Delta + s/6$ has trivial kernel;  the latter is
automatic if  the conformal class $[g]$ is contains a metric with $s> 0$, 
never happens if $[g]$ is contains a metric with
$s\equiv 0$, and may or may not happen if $[g]$  contains a metric with $s<0$. 
Since the operator $\Delta + s/6$ is self-adjoint, it also has trivial 
cokernel, and the equation 
$$(\Delta + s/6) u =f$$
therefore has a unique smooth solution $u$ for any smooth function $f$;
it follows that it also has a unique distributional solution $u$ for any distribution $f$. 
If $y\in M$ is any point, and if $\delta_y$ is the Dirac delta distribution 
centered at $y$, we thus have a unique distributional solution $G_y$
of the equation 
$$(\Delta + s/6)G_y =\delta_y.$$
Since $\delta_y$ is  identically zero on $M-\{ y\}$, elliptic
regularity tells us that $G_y$ is actually smooth away from $y$. 
In general, one has an expansion
$$
G_y =  \frac{1}{4\pi^2}\frac{1}{\varrho^2}+ O(\log \varrho )
$$
near  $\varrho$ denotes the distance from $y$,
but when  $(M,g)$ is anti-self-dual it in fact turns out  \cite{atgrn} that 
$$G_y =  \frac{1}{4\pi^2}\frac{1}{\varrho^2}+ \mbox{ bounded terms. }$$
In this article, the function $G_y$ will be called the {\em conformal Green's function} of $(M,g,y)$.

The motivation for this terminology is that the 
 Yamabe Laplacian is  a {\em conformally invariant}
differential operator when viewed as a map between sections of suitable real line bundles; 
the geometric reason for this is that for any smooth
function $u\not\equiv 0$, the expression $6u^{-3}(\Delta + s/6)u$
computes the scalar curvature of the conformally related metric $u^2g$ on the open set
$u\neq 0$. One useful consequence of this is that, for any smooth  function
$u > 0$ on $M$, a  constant times $u^{-1}G_y$ 
is the conformal Green's function of $(M,u^2g,y)$. 

Now the celebrated proof of the Yamabe conjecture \cite{lp}
tells us that any conformal class on any compact manifold contains metrics of constant 
scalar curvature. In particular, any conformal class contains 
metrics whose scalar curvature has the same sign at every point. But actually, this last
assertion is much more elementary. Indeed, if $u\not\equiv 0$ is an eigenfunction corresponding 
to the lowest eigenvalue $\lambda$ of the Yamabe Laplacian, then $u$ has empty nodal 
set, and $u^2g$ is therefore a conformally related metric whose scalar curvature
has the same sign as $\lambda$ everywhere on $M$. 
Similar considerations also show that if two metrics with scalar curvatures
of fixed signs are conformally related, then their scalar curvatures have
the {same} sign. 
The {\em sign of Yamabe constant}
of a conformal class, meaning the sign of the constant scalar curvature of the metric produced by the 
proof of the Yamabe conjecture, therefore coincides with the sign of the smallest Yamabe
eigenvalue $\lambda$
 for any metric in the conformal class.

 Here is another way of determining the sign of the Yamabe constant:  
 
 \begin{lem}\label{greenjeans}
 Let $(M,g)$ be a compact Riemannian $4$-manifold
 whose Yamabe Laplacian $\Delta + s/6$ has trivial kernel.
Let $y\in M$ be any point. 
Then the conformal class $[g]$ contains a metric 
of positive scalar  curvature if and only if $G_y(x) \neq 0$
for all $x\in M - \{ y\}$. Moreover, if $[g]$ contains a metric 
of negative scalar  curvature, then $G_y(x) < 0$ for some
$x\in M$.  

These assertions also hold for any finite sum $G_{y_1}+\cdots + G_{y_m}$,
or for any other finite linear combinations of conformal Green's functions with 
positive coefficients. 
\end{lem}
 \begin{proof}
 Since 
 the Yamabe Laplacian is conformally invariant when viewed as 
 acting on functions of the appropriate conformal weight,
 we may assume from the outset that either $s>0$ everywhere, or else
$s < 0$  everywhere. 
 
Now notice that notice that
 $$\frac{1}{6}\int sG_y~d\mu = \int   (\Delta + s/6) G_y~d\mu = \int \delta_y ~d\mu = 1 > 0.$$ 
 Thus, if $s < 0$, $G_y$ must be negative somewhere, and since 
 $G_y\to +\infty$ at 
 $y$, the Green's function 
  must also have 
 a zero by continuity. 
 Notice that the same argument works for any finite linear combinations of Green's functions
 with positive coefficients.

 On the other hand, $G_y^{-1}((-\infty , a ] )$ is compact for
 any $a \in \RR$, and it follows that  $G_y$ has a minimum. 
 But if $s > 0$, then   $G_y = \frac{1}{s} \nabla \cdot \nabla G_y$ on $M -\{ y\}$,
 and hence $G_y \geq 0$ at its minimum. Moreover, if the minimum were actually
 zero, we could apply  Hopf's strong maximum principle \cite{giltrud,protter} to $-G_y$, and conclude 
 that $G_y\equiv 0$,  contradicting the fact that $G_y\to \infty$ at 
 $y$. Thus $G_y > 0$ everywhere, as claimed.  As a consequence, any  finite linear combinations of Green's functions
 with positive coefficients is a sum of positive functions, and so is positive at all points
 where its value is defined. 
  \end{proof}

Now suppose that $(M,g)$ is a compact {\em anti-self-dual} Riemannian
$4$-manifold, and let $Z$ be its twistor space.  If $U\subset M$ is any open 
subset, and if $Z_U\subset Z$ is its inverse image in the twistor space,
the Penrose transform  \cite{bailsing,hitka}  gives a natural one-to-one correspondence 
between  $H^1 (Z_U, {\mathcal O}(K^{1/2}))$ and the smooth 
complex-valued functions
on $U$ which solve $(\Delta  + s/6)u=0$. 
Given a cohomology class $\psi \in H^1 (Z_U, {\mathcal O}(K^{1/2}))$, the value of 
the corresponding function $u_\psi$ at $x\in U$ is obtain by restricting $u_\psi$
to the real twistor line $P_x\subset Z$ to obtain an element of
$H^1(P_x, {\mathcal O}(K_Z^{1/2}))\cong H^1(\bcp_1 , {\mathcal O}(-2) ) \cong \CC$. 
Note that  $u_\psi$ is  ostensibly only a section of a line bundle, but the choice
of a metric $g$ in the conformal class turns out to determine a canonical  trivialization
of  this line
bundle \cite{hitka}, and $u_\psi$ then becomes a  function in the ordinary sense.

In particular, our compact anti-self-dual $4$-manifold $(M,g)$ satisfies
$\ker (\Delta + s/6)=0$
iff its twistor space $Z$ satisfies $H^1(Z, {\mathcal O}(K^{1/2}))=0$, 
and Serre duality tells us that the latter happens
iff $H^2(Z, {\mathcal O}(K^{1/2}))=0$.
When any of these thee equivalent conditions is met, we
then have a conformal 
 Green's function 
$G_y$ for any chosen $y\in M$, and $G_y$ then 
corresponds to a particular element of $H^1(Z-P_y, {\mathcal O}(K^{1/2}))$,
where $P_y$ is the twistor line corresponding to  $y$. 
What is this  mysterious cohomology class?
The answer was discovered by  Atiyah \cite{atgrn}, and involves
 a construction largely due to Serre \cite{serremod} and Horrocks \cite{horrocks}: 

\begin{lem} \label{sehr}
Let $W$ be a (possibly non-compact) complex manifold, and let $V\subset W$ be a 
closed complex submanifold of complex codimension $2$. Let $N\to V$
denote the normal bundle $TW/TV$ of $V$, and suppose that there is 
a holomorphic line bundle $L\to W$ such that 
\begin{itemize}
\item 
$L|_V\cong \wedge^2 N$;
\item $H^1(W, {\mathcal O}(L))=0$; and 
\item $H^2(W, {\mathcal O}(L))=0$.
\end{itemize}
Then there is a rank-$2$ holomorphic vector bundle $E\to W$,
together with a holomorphic section $\zeta \in \Gamma (W, {\mathcal O}(E))$
such that 
\begin{itemize}
\item 
$\wedge^2 E \cong L$;
\item $\zeta =0$ exactly at $W$; and 
\item $d\zeta : N\to E$ is an isomorphism.
\end{itemize}
This $(E, \zeta)$ is unique up to isomorphism if we also demand that 
the isomorphism $\det d\zeta: \wedge^2 N\to \wedge^2E|_V$ should agree with a given
 isomorphism $ \wedge^2 N\to L|_V$. 
 The pair  $(E, \zeta)$
gives rise to an extension 
$$0\to {\mathcal O}(L^*)\to {\mathcal O}(E^*)
\stackrel{\cdot \zeta}{\longrightarrow} {\mathscr I}_V\to 0,$$
where ${\mathscr I}_V$ is the ideal sheaf of $V$; and by restriction to 
$W-V$, this extension  determines 
an element of $H^1(W-V,  {\mathcal O}(L^*))$. 
\end{lem}
\begin{proof}
Because \cite{altkle,GH} $V\subset W$ is smooth and of codimension $2$, 
$${\mathcal Ext}^q({\mathcal O}_V, {\mathcal O}(L^*))=
\begin{cases}
      0& q=0,1, \\
      {\mathcal O}_V(L^*\otimes \wedge^2N)& q=2,
\end{cases}  
      $$
and the spectral sequence
$$E_2^{p,q}= H^p(W, {\mathcal Ext}^q({\mathcal O}_V, {\mathcal O}(L^*)))
\Longrightarrow \text{\bf Ext}^{p+q}_W ({\mathcal O}_V, {\mathcal O}(L^*))$$
therefore tells us that 
$$\text{\bf Ext}^{2}_W ({\mathcal O}_V, {\mathcal O}(L^*))= \Gamma (V,  
{\mathcal O}_V(L^*\otimes \wedge^2N)).$$
 On the other hand, the tautological short exact sequence of sheaves
$$
0\to {\mathscr I}_V\to {\mathcal O}\to {\mathcal O}_V\to 0
$$
induces a long exact sequence

\begin{eqnarray*}
 \cdots \to& \text{\bf Ext}^{1}_W ({\mathcal O}, {\mathcal O}(L^*)) & 
 \to  \text{\bf Ext}^{1}_W ({\mathscr I}_V, {\mathcal O}(L^*))
 \\
\to \text{\bf Ext}^{2}_W ({\mathcal O}_V, {\mathcal O}(L^*))
\to &  \text{\bf Ext}^{2}_W ({\mathcal O}, {\mathcal O}(L^*))
 & \to \cdots
\end{eqnarray*}
and since 
$$\text{\bf Ext}^{j}_W ({\mathcal O}, {\mathcal O}(L^*))
=H^j(W, {\mathcal O}(L^*))$$
 is  assumed to vanish when $j=1,2$, the Bockstein map of this 
long exact sequence therefore gives us an isomorphism 
$$\text{\bf Ext}^{1}_W ({\mathscr I}_V, {\mathcal O}(L^*))
\cong \Gamma (V,  
{\mathcal O}_V(L^*\otimes \wedge^2N)).$$
In particular,  any choice of isomorphism  $\wedge^2N\to L|_V$  gives us an 
extension 
\begin{equation}
\label{serre}
0\to {\mathcal O}(L^*) \to {\mathcal O}(E^*)\to {\mathscr I}_V\to 0
\end{equation}
of sheaves on $W$. 
The class of this extension is called the 
{\em Serre class} $\lambda (V) \in \text{\bf Ext}^{1}_W ({\mathscr I}_V, {\mathcal O}(L^*))$,
and its restriction to $W-V$ is an element of $H^1(W-V,  {\mathcal O}(L^*))$. 
Strictly speaking, the
 Serre class  depends on a choice of isomorphism $\wedge^2N\to L|_V$,
but any two such extensions are intertwined by an automorphism of $L|_V$. 

Since  any   isomorphism $\wedge^2N\to L|_V$ corresponds to a   section 
of ${\mathcal Ext}^2({\mathcal O}_V, {\mathcal O}(L^*))$ which is 
non-zero at each point of $V$, the corresponding  extension (\ref{serre}) 
is locally free, with $\Lambda^2E^* = L^*$. Tensoring the inclusion ${\mathcal O}(L^*) \hookrightarrow {\mathcal O}(E^*)$ by  $L$, 
we thus obtain an inclusion $${\mathcal O} \hookrightarrow {\mathcal O}(E^*\otimes
\Lambda^2 E)= {\mathcal O}(E),$$ and the image of $1$ under this map
is then a section $\zeta \in \Gamma (W, {\mathcal O}(E))$
with all the advertised properties. 
\end{proof}

\begin{prop}[Atiyah] \label{atiyah}
Let $(M,g)$ be a compact anti-self-dual $4$-manifold 
with twistor space $Z$, and assume that $(M,g)$ has 
 $$\ker (\Delta + s/6)=0.$$ Let $y\in M$ be any point, and
 let $P_y\subset Z$ be the corresponding twistor line. 
 Then the image of the Serre class $\lambda (P_y)\in 
 \text{\bf Ext}^{1}_Z ({\mathscr I}_{P_y}, {\mathcal O}(K^{1/2}))$
 in 
 $H^1(Z-P_y, {\mathcal O}(K^{1/2}))$ is the Penrose transform of the Green's function $G_y$
 times a non-zero constant. 
\end{prop}

Indeed, if one identifies $K^{1/2}_Z|_{P_y}$ with $K_{P_y}$ according the isomorphism
determined by $g$ and the conventions of \cite{hitka},  the relevant constant 
turns out to be exactly $4\pi$. 

Combining this remarkable result with Lemma \ref{greenjeans} now gives us a twistorial
criterion for determining  whether an  anti-self-dual conformal class 
has positive Yamabe constant: 

\begin{prop}\label{nice}
Let $Z$ be the twistor space of a  compact anti-self-dual $4$-manifold $(M, [g])$,
and let $P_y\subset Z$ be a real twistor line.  Then the conformal class
$[g]$ contains a metric $g$ of positive scalar curvature if and only if
\begin{itemize}
\item $H^1(Z, {\mathcal O}(K^{1/2}))=0$,  and 
\item the 
holomorphic vector bundle $E\to Z$ with $\wedge^2E\cong K^{-1/2}$ 
associated to $P_y$  by   Lemma \ref{sehr} satisfies 
$E|_{P_x}\cong {\mathcal O}(1)\oplus {\mathcal O}(1)$
for every real twistor line $P_x$.
\end{itemize}
\end{prop}
\begin{proof} Let us first recall that 
$$
\ker (\Delta + s/6) = H^1(Z,  {\mathcal O}(K^{1/2})) = [H^2(Z,  {\mathcal O}(K^{1/2}))]^*,$$
so that a necessary condition for the positivity of the Yamabe constant is certainly
the vanishing of $H^1(Z, {\mathcal O}(K^{1/2}))$. When this happens, 
Lemma \ref{sehr} then allows us to construct $E\to Z$. 
On $Z-P_y$, $E$ is then given by an extension 
$$0\to {\mathcal O}\to  {\mathcal O}(E)\to {\mathcal O}(K^{-1/2})\to 0,$$
and this extension is represented by an element of $H^1(Z-P_y,  {\mathcal O}(K^{1/2}))$.
The value of the  Penrose transform of this class  at $x\neq y$ is obtained via
the restriction map 
$$
H^1(Z-P_y,  {\mathcal O}(K^{1/2}))\to H^1(P_x,  {\mathcal O}_{P_x}(K^{1/2}))\cong \CC
$$
and its value at $x$ is therefore non-zero iff the induced extension 
$$0\to {\mathcal O}\to E|_{P_x}\to {\mathcal O}(2)\to 0$$
does not split; and this  happens iff $E|_{P_x}\cong {\mathcal O}(1)\oplus {\mathcal O}(1)$. 
Since we also have $E|_{P_y}\cong {\mathcal O}(1)\oplus {\mathcal O}(1)$
by construction, 
the result now follows from 
lemma \ref{greenjeans} and Proposition \ref{atiyah}. 
\end{proof}


\section{The Sign of the Scalar Curvature}
\label{scalar}

We are now ready to approach  the problem of determining  the sign of the 
Yamabe constant for the anti-self-dual conformal classes constructed in Theorem \ref{thmb}.
The results obtained in this section  are loosely 
inspired by the  work of  Dominic Joyce   \cite{joyscal}
on the Yamabe constants of certain conformal classes on connect sums,
although the techniques employed here are  completely different from Joyce's.

\begin{lem}
Let $\varpi : {\mathcal Z} \to {\mathcal U}$ be a $1$-parameter standard deformation of 
$Z_0$, where $Z_0$ is as  in Theorem \ref{thmb}, and ${\mathcal U}\subset \CC$ is an open disk
about 
the origin. 
Let $k\tilde{Z}_{EH}$ be the union of the Eguchi-Hanson components of 
$Z_-\subset 
Z_0$, which is a 
 non-singular
complex  hypersurface in $\mathcal Z$,
and let ${\mathscr I}_{k\tilde{Z}_{EH}}\subset {\mathcal O}$
denotes its ideal sheaf.
 Then the invertible sheaf 
${\mathscr I}_{k\tilde{Z}_{EH}}(K_{\mathcal Z})\subset {\mathcal O}_{\mathcal Z}(K_{\mathcal Z})$ has
a square-root as a holomorphic line bundle.
\end{lem}
\begin{proof}
A holomorphic line bundle has a holomorphic square-root iff its second Stieffel-Whitney class
$w_2\in H^2(\ZZ_2)$ vanishes. Write $\mathcal Z$ as $U\cup V$, where $U$ is 
a tubular neighborhood of $Z_+$, $V$ is a tubular neighborhood of $Z_-$, and 
$U\cap V$ is a tubular neighborhood of $Q$, so that these open sets deform retract to 
$Z_+$, $Z_-$, and $Q$, respectively. Since each component of $Q$ is simply connected, 
the Mayer-Vietoris sequence
$$
\cdots \to H^1(U\cap V, \ZZ_2) \to H^2(U\cup V, \ZZ_2 ) \to H^2 (U, \ZZ_2) 
\oplus H^2 (V, \ZZ_2)\to
\cdots 
$$
therefore tells us that it is enough to check that the restrictions of our line bundle 
to $Z_+$ and $Z_-$ both have square-roots.

It thus suffices to produce an explicit square-root of the restrictions 
of ${\mathscr I}_{k\tilde{Z}_{EH}}\otimes K_{\mathcal Z}$
to each copy of $\tilde{Z}_{FS}$, each copy of $\tilde{Z}_{EH}$,
and to $Z_+$. 
On each copy of $\tilde{Z}_{FS}$, such a square-root  is given by $[Q]\otimes K^{1/2}_{Z_{FS}}$,
where $[Q]$ is the divisor of the exceptional quadric, and where 
$Z_{FS}^{1/2}$ is the  pull-back of $K^{1/2}$ from the twistor space $Z_{FS}$ via
the blowing-down map. On each copy of $\tilde{Z}_{EH}$, such a square-root  is given by 
the pull-back of ${\mathcal O}(-2)$ via the projection $\tilde{Z}_{EH}\to \bcp_1$.
And on $Z_+$, there is a natural choice of square-root whose sections
are the $\ZZ_2$ invariant sections of $K^{1/2}_Z$, pulled-back to the blow-up
$\tilde{Z}$, twisted by the divisors $Q_{k+1}, \ldots , Q_{k+\ell}$. 
That each of these bundles really has the correct square
 can be verified directly  using the adjunction formula; the 
details are left as an exercise for the interested reader. 
\end{proof}

\begin{lem}\label{family}
Let $\varpi : {\mathcal Z} \to {\mathcal U}$ be a $1$-parameter standard deformation of 
$Z_0$, where $Z_0$ is as in Theorem \ref{thmb}, and ${\mathcal U}\subset \CC$ is a neighborhood
of the origin. Let $L\to {\mathcal Z}$ be the holomorphic line bundle defined by
$${\mathcal O}(L^*)
= [{\mathscr I}_{k\tilde{Z}_{EH}}(K_{\mathcal Z})]^{1/2}\otimes {\mathscr I}_{\ell\tilde{Z}_{FS}},$$
where the hypersurface 
$\ell\tilde{Z}_{FS}\subset {\mathcal Z}$ is the union of the Fubini-Study 
components of $Z_-$. 
If the twistor space $Z$ of $(Y,[g])$ satisfies $H^1(Z, {\mathcal O}(K^{1/2}))=0$,
then  by possibly replacing  ${\mathcal U}$ with a smaller 
 neighborhood  of $0\in \CC$ and simultaneously replacing 
 ${\mathcal Z} $ with its inverse image, we can arrange for
 our 
  complex $4$-fold $\mathcal Z$ 
to satisfy 
$$H^1 ({\mathcal Z} , {\mathcal O}(L^*))=H^2 ({\mathcal Z} , {\mathcal O}(L^*))=0.$$
\end{lem}
\begin{proof} Since any open set in $\CC$ is Stein, the Leray spectral sequence 
tells us that it would suffice to show that the direct image sheaves $\varpi_*^j{\mathcal O}(L^*)$
vanish for $j=1,2$. But since $\varpi$ is flat and we are allowed to shrink ${\mathcal U}$ if 
necessary, semi-continuity \cite{bast} asserts that it is enough to show that 
$H^j(Z_0,  {\mathcal O}(L^*))=0$ for $j=1,2$. 

The normalization of $Z_0$ is the disjoint union $Z_+\sqcup Z_-$, and
we have an exact sequence
$$0\to {\mathcal O}_{Z_0}(L^*)\to \nu_*{\mathcal O}_{Z_+}(L^*)\oplus
\nu_*{\mathcal O}_{Z_-}(L^*)
\to {\mathcal O}_Q(L^*)\to 0$$
where $\nu : Z_+\sqcup Z_-\to Z_0$ is the identification map. However,
$ {\mathcal O}_{Z_0}(L^*)$ exactly consists of $\ZZ_2$-invariant sections
of the pull-back of $K^{1/2}_Z$, and the Leray spectral sequence 
therefore tells us that $$H^j(Z, {\mathcal O}(K^{1/2}))=0\Longrightarrow 
H^j(Z_+,  {\mathcal O}(L^*))=0,$$
so our vanishing hypothesis guarantees that these groups vanish for all $j$. 
On the other hand, 
$$
H^j(\tilde{Z}_{EH}, {\mathcal O}(L^*))
= H^j(\tilde{Z}_{FS}, {\mathcal O}(L^*))
= H^j(Q_j, {\mathcal O}(L^*))
=\begin{cases}
   \CC   &j=1, \\
     0 & j\neq 1.
\end{cases}
$$
and each of the relevant restriction maps, from $H^1(\tilde{Z}_{EH}, {\mathcal O}(L^*))$ or 
$H^1(\tilde{Z}_{FS}, {\mathcal O}(L^*))$ to the cohomology group 
$H^1(\bcp_1\times \bcp_1 , {\mathcal O}(-2,0))$ of the appropriate
quadric $Q_j$, is an isomorphism. 
Hence $H^j(Z_0, {\mathcal O}(L^*))=0$ for all $j$, and the result follows. 
\end{proof}

Now choose a real twistor line $P_{x}\in Z_+$, and extend this as a $1$-parameter family 
of twistor lines  in $P_{x_{\mathfrak u}}\in Z_{\mathfrak u}$ for ${\mathfrak u}$ 
near $0\in \CC$ and such that 
 $P_{x_{\mathfrak u}}$ is a real twistor line 
for ${\mathfrak u}$ real. By possibly shrinking ${\mathcal U}$, we may then arrange that 
${\mathcal P}= \cup_{\mathfrak u} P_{x_{\mathfrak u}}$ is a
 closed submanifold of ${\mathcal Z}$ and that 
 $H^1 ({\mathcal Z} , {\mathcal O}(L^*))=H^2 ({\mathcal Z} , {\mathcal O}(L^*))=0$.
 The hypotheses of Lemma \ref{sehr} are then satisfied, and we thus obtain 
 a holomorphic vector bundle $E\to {\mathcal Z}$  and a holomorphic section
 $\zeta$ vanishing exactly along $\mathcal P$; moreover, the corresponding extension
 $$0\to {\mathcal O}(L^*)\to {\mathcal O}(E^*)\to {\mathscr I}_{\mathcal P}\to 0$$
 gives us an element of $\lambda({\mathcal P})\in H^1({\mathcal Z}-{\mathcal P}, {\mathcal O}(L^*))$.
  Since the restriction of 
 $L^*$ to any smooth fiber $Z_{\mathfrak u}$, ${\mathfrak u}\neq 0$,
  is just $K^{1/2}$, Proposition \ref{atiyah} tells us that 
  the restriction of $\lambda ({\mathcal P})$ to 
 $Z_{\mathfrak u}$, ${\mathfrak u} > 0$, has Penrose transform equal to   a positive 
 constant times the conformal Green's function
 of $(\tilde{X}\# \ell \overline{\bcp}_2, g_{\mathfrak u}, x_{\mathfrak u})$ 
 for any ${\mathfrak u} > 0$. However,  we may also restrict $(E,\zeta )$
  to  $Z_+$, and, by pulling-back and pushing down, convert this into
  a $\ZZ_2$-invariant holomorphic vector bundle on the twistor space 
  $Z$ of $(Y,g)$. This bundle on $Z$ then has determinant line bundle $K^{-1/2}$, and 
  comes equipped with a section 
  vanishing exactly at the twistor lines of the two pre-images of $y_1, y_2$
  of $x\in X$; by Proposition \ref{atiyah},  the Penrose transform of this object corresponds, 
  according to your taste, either to  the Green's function $G_x$
 on$(X,g)$ or to the sum $G_{y_1}+G_{y_2}$ on $(Y,g)$. 
 If $g$ has negative scalar curvature, Lemma \ref{greenjeans} 
 thus tells us that 
 is a region of
 $X$ where $G_x < 0$, and defroming the twistor lines of this into  
 $Z_{\mathfrak u}$ for small  ${\mathfrak u}> 0$ then shows that the conformal Green's function 
 of $[g_{\mathfrak u}]$ is  negative somewhere for any small $\mathfrak u$.
 By  Lemma \ref{greenjeans}, we thus obtain the following:

\begin{thm}\label{neg}
In  Theorem \ref{thmb}, suppose that 
 $(Y,g)$ is an anti-self-dual manifold with  $s<0$
and $\ker (\Delta + s/6)=0$. 
Then for all sufficiently small ${\mathfrak u}>0$, the conformal class  $[g_{\mathfrak u}]$ contains
 a metric of   negative scalar curvature.
\end{thm}

The positive case is similar, but is slightly more delicate. Instead of just restricting
$\lambda ({\mathcal P})$ on rational curves in $Z_+$, we must also consider
what happens when we restrict this class to twistor lines in $Z_-$. However, 
we already saw in the proof of Lemma \ref{family}, an element of $H^1(\tilde{Z}_{EH},
{\mathcal O}(L^*))\cong \CC$ or $H^1(\tilde{Z}_{FS},
{\mathcal O}(L^*))\cong \CC$ is non-zero iff its restriction to the corresponding exceptional
quadric is non-zero, and this has the effect that the restriction of the 
cohomology class to every twistor line in either of these spaces  is non-zero if
there is a rational curve in the quadric on which the class is non-zero. Thus, 
when the conformal Green's function $G_x$ of $(X,g)$ is positive, 
the vector bundle $E$ determined by 
$\lambda ({\mathcal P})$ has splitting type ${\mathcal O}(1)\oplus {\mathcal O}(1)$ 
on all  the $\sigma_0$-invariant rational curves in $Z_0$
which are limits of real twistor lines in ${\mathcal Z}_{\mathfrak u}$ as ${\mathfrak u}\to 0$.
It therefore has the same splitting type on all the real twistor lines of ${\mathcal Z}_{\mathfrak u}$
for ${\mathfrak u}$ small, and Proposition \ref{nice}  therefore tells us:

\begin{thm}\label{pos}
In  Theorem \ref{thmb}, suppose 
that $(Y,g)$ is an anti-self-dual manifold with  $s>0$.
Then for all sufficiently small ${\mathfrak u}>0$, the conformal class  $[g_{\mathfrak u}]$ contains
 a metric of  positive scalar curvature.
\end{thm}

In this positive case, it is  interesting to re-examine the above construction 
in purely Riemannian terms. In this setting, the positivity of the Green's functions allows us
to define a family of asymptotically flat, scalar-flat, anti-self-dual metrics 
 $\tilde{g}_{\mathfrak u}=  G_{x_{\mathfrak u}}^2g_{\mathfrak u}$ on 
$(\tilde{X}\# \ell \overline{\bcp}_2) - \{ pt\}$.
What the above construction tells us is that these metrics  converge,
in the pointed Gromov-Hausdorff sense \cite{gromhaus}, 
 to the orbifold  metric $\tilde{g}=G_x^2g$ on $X-x$. However, 
there is something else going on in certain regions, where  viewing these metrics
under higher and higher magnification results in  a family that  
 converges to the
Eguchi-Hanson metric or to the {\em Burns metric}, meaning the Green's
function rescaling of the Fubini-Study metric on $\overline{\bcp}_2-\{pt\}$. 
It is the appearance of the ideal sheaves in the definition of $L^*$ which
accounts for the fact that these seemingly incompatible pictures  simultaneously  apply
at wildly different length-scales.

\section{Conformally Flat Orbifolds}
\label{limitset}

Consider the involution of  $S^1 \times S^3\subset \CC \times \HH$ given
by $(z,q)\mapsto (\bar{z},\bar{q})$. 
This involution has only four fixed points, namely $(z,q)=(\pm 1 , \pm 1)$; 
 and near each of these isolated fixed points, the
involution necessarily  looks exactly like reflection through the origin in $\RR^4$. 
We can therefore construct an involution of 
the connected sum 
$(S^1 \times S^3)\# (S^1 \times S^3)$
by cutting out a ball centered at a fixed point of the involution of
each copy of $S^1 \times S^3$, and then being careful to carry out the usual  gluing procedure 
in a  $\ZZ_2$-equivariant manner. The resulting involution 
$$\phi : (S^1 \times S^3)\# (S^1 \times S^3)\to (S^1 \times S^3)\# (S^1 \times S^3)$$
then has exactly $6$ fixed points, and
may usefully be thought of as a  $4$-dimensional analog of the hyperelliptic involution 
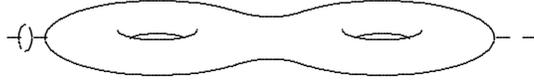
\begin{center}
\mbox{
\beginpicture
\setplotarea x from 40 to 290, y from -5 to 65
\ellipticalarc axes ratio 3:1  270 degrees from 150 40
center at 120 30
\ellipticalarc axes ratio 3:1  -270 degrees from 175 40
center at 205 30
\ellipticalarc axes ratio 4:1 -180 degrees from 135 33
center at 120 33
\ellipticalarc axes ratio 4:1 145 degrees from 130 30
center at 120 29
\ellipticalarc axes ratio 4:1 180 degrees from 190 33
center at 205 33
\ellipticalarc axes ratio 4:1 -145 degrees from 195 30
center at 205 29
\ellipticalarc axes ratio 1:2 140 degrees from 70 35
center at 70 30
\ellipticalarc axes ratio 1:2 140 degrees from 70 25
center at 70 30
{\setquadratic 
\plot 150 40    163 38    175 40   /
\plot 150 20    163 22    175 20  /
}
{\setlinear
\setdashes 
\plot  60 30 80 30   /
\plot  245 30 265 30   /
}
\endpicture
}
\end{center}
of a Riemann surface of genus $2$.  

In this section, we will be interested
in conformally flat orbifold metrics on $X=[(S^1 \times S^3)\# (S^1 \times S^3)]/\ZZ_2$,
or in other words, in $\phi$-invariant, conformally flat metrics on $(S^1 \times S^3)\# (S^1 \times S^3)$.
The  key result we'll need
is the following:  

\begin{prop}\label{cle}
There is a real-analytic  family $g_t$, $t\in (0,1)$, of  Riemannian 
metrics on $(S^1 \times S^3)\# (S^1 \times S^3)$ with the following properties:
\begin{itemize}
\item for each $t$, the metric $g_t$ is locally conformally flat; 
\item  for each $t$, the involution $\phi$ is an isometry of  $g_t$; 
\item for each $t$, the scalar curvature $s$ of $g_t$ has a fixed sign; 
\item when $t$ is sufficiently close to $0$, $g_t$ has $s > 0$; 
\item when $t$ is sufficiently close to $1$, $g_t$ has $s < 0$; 
\item the set of $t$ for which $\ker (\Delta + s/6)\neq \{ 0\}$ 
is discrete; and 
\item there are only finitely many values of $t$ for which $g_t$ has $s\equiv 0$.
\end{itemize}
Consequently, there is a  $t_0\in (0,1)$ and an $\varepsilon > 0$ 
such that, for all  $t\in (t_0-\varepsilon , t_0+ \varepsilon )$,
 the scalar curvature of  $g_t$ has the same sign as $t_0 -t$,
 and such that the Yamabe Laplacian $\Delta + s/6$ has trivial
 kernel for all $t\in (t_0-\varepsilon , t_0)\cup (t_0 , t_0+ \varepsilon )$. 
 \end{prop}

Our proof of the existence of  such a family hinges on a result of 
Schoen and Yau \cite{schyaudim}, and the  construction used here 
is analogous to related constructions of   Kim  \cite{jongsu} and Nayatani \cite{shinpat}. 
Let us begin by observing that 
 $(S^1 \times S^3)\# (S^1 \times S^3)$ can be obtained from 
 $S^4$ by deleting four  balls, and identifying 
 the resulting boundary spheres in pairs via reflections. 
Now think of $S^4$ as 
 $\HH\PP_1=\HH \cup \{ \infty \}$, and 
 let $D_t\subset \HH \cup \{ \infty \}$
 be the complement of the four open balls
 $$B_t(\sqrt{3}), B_t(-\sqrt{3}), B_t(i) , B_t(-i)  \subset  \HH$$
  of radius $t$;  and we henceforth  
  stipulate that  $t < 1$,  so as to guarantee that the closures of these four balls are pairwise disjoint. 
  We may then think of  $(S^1 \times S^3)\# (S^1 \times S^3)$ 
  as obtained from $D_t$ by 
 identifying $\partial B_t(i)$ with $\partial B_t(-i)$
 via the reflection $q\mapsto \bar{q}$, and  identifying 
$\partial  B_t(\sqrt{3})$ with $\partial  B_t(-\sqrt{3})$
  via the reflection $q\mapsto -\bar{q}$:
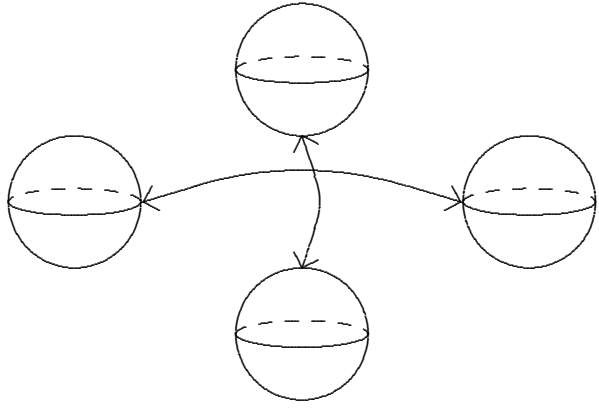
  \begin{center}
\mbox{
\beginpicture
\setplotarea x from -150 to 150, y from -100 to 100
 \circulararc 360 degrees from 61 0 center at  86 0
  \circulararc 360 degrees from -61 0 center at  -86 0
\circulararc 360 degrees from 0 25 center at  0 50
\circulararc 360 degrees from 0 -25 center at  0 -50
\ellipticalarc axes ratio 5:1 180 degrees from 61 0 center at 86 0
\ellipticalarc axes ratio 5:1 -180 degrees from -61 0 center at -86 0
\ellipticalarc axes ratio 5:1 180 degrees from -25 50 center at 0 50
\ellipticalarc axes ratio 5:1 180 degrees from -25 -50 center at 0 -50
\arrow <5pt> [1,2] from 40 6   to 60 0 
\arrow <5pt> [1,2] from 5 -10   to 0 -25
\arrow <5pt> [1,2] from -40 6   to -60 0 
\arrow <5pt> [1,2] from 5 10   to 0 25
{\setdashes
\ellipticalarc axes ratio 5:1 -180 degrees from 61 0 center at 86 0
\ellipticalarc axes ratio 5:1 180 degrees from -61 0 center at -86 0
\ellipticalarc axes ratio 5:1 -180 degrees from -25 50 center at 0 50
\ellipticalarc axes ratio 5:1 -180 degrees from -25 -50 center at 0 -50
}
{
\setquadratic
\plot  -47 6  -7 12    33 6 /
\plot  -2 -10  0   0   -2 10 /
}
\endpicture
}
\end{center}
 Since
the reflections that we have used to identify
the boundary components in pairs are just the restrictions to the relevant spheres 
of the global conformal transformations 
$$q\mapsto t^2 (q+i)^{-1} + i  ~~\mbox{and}~~q\mapsto -t^2 (q -\sqrt{3})^{-1}-\sqrt{3}$$
 of $S^4= \HH \cup \{ \infty\}$, we  may therefore define
 a  unique flat conformal metric $[g_t]$ on $(S^1 \times S^3)\# (S^1 \times S^3)$
 by restricting the standard conformal metric on $S^4$ to $D_t$, and then
 pushing this structure down to  $D_t/\sim$. 
 
To complete the picture, we let $\phi$  act on 
$D_t/\sim$ by 
$q\mapsto -q$. Note that we obviously have  $\phi^*[g_t] = [g_t]$. 
The six fixed points of $\phi$  are just  
$ i\pm t, \sqrt{3} \pm t,0$  and $\infty$.
That this  $\phi$ coincides with the previously-described 
involution of $(S^1 \times S^3)\# (S^1 \times S^3)$
 may be seen, when
 $t \in (0,1/3)$,  by first decomposing $\HH \cup \{\infty \}$ 
 into the hemispheres $\| q \| \leq 4/3$ and  $\| q \| \geq 4/3$. 
 Thus $\phi$ can be constructed by first letting $q\mapsto -q$ act on 
 two separate copies of $S^4$ minus  {\em two} balls,  each with its boundary
 components identified via  reflections, and then forming the $\ZZ_2$-equivariant 
 connected sum of these manifolds; but each of these two building blocks looks like  
$ S^3\times\left([-1,1]/\{-1, 1\}\right)$  equipped with the involution
$(q,t)\mapsto (\bar{q},-t)$, so the claim follows.

So far, we have only constructed a  family $[g_t]$ of $\phi$-invariant flat conformal classes,
but we next need to worry about how nicely these conformal structures vary with $t$. 
However, it is not hard to see that they are real-analytic in $t$, since on a given open 
neighborhood $U$ of a given $D_{t^\prime}\subset S^4$, we are are simply gluing together 
 neighborhoods of the boundary spheres via the the M\"obius transformations
$$q\mapsto t^2 (q+i)^{-1} +i , ~~\mbox{and}~~q\mapsto -t^2 (q -\sqrt{3})^{-1}-\sqrt{3}$$
for $t$ near $t^\prime$, and these transformations 
depend real-analytically on $t$. Since the sheaf of real-analytic functions 
is acyclic \cite{gunros}, we can now choose a real-analytic family of $\phi$-invariant 
metrics $h_t$ which represents the family of conformal classes $[g_t]$. 
Let $\lambda_t$ be the smallest eigenvalue of the Yamabe Laplacian $\Delta_{h_t}+ s_{h_t}/6$
of $h_t$, and let $f_t$ be an eigenfunction of eigenvalue $\lambda_t$ and integral $1$.
By the minimum principle \cite{giltrud}, $f_t$ is everywhere positive, and it
follows that it must be unique; in particular, $f_t$ must be $\phi$-invariant. 
Moreover, this uniqueness tells us that  $\lambda_t$ has multiplicity $1$.
Hence $\lambda_t$ never meets another eigenvalue as $t$ varies, so 
perturbation theory \cite{kato} tells us that $\lambda_t$ and  $f_t$ depend
real-analytically on $t$. 
Now set 
$$g_t = f_t^2 h_t,$$
and notice that the scalar curvature 
$$s_{g_t}= f_t^{-3}(6\Delta_{h_t}+ s_{h_t}) f_t = 6\lambda_t f_t^{-2}$$
of this metric 
has the same sign as $\lambda_t$ at every point. 
Thus $g_t$ is a real-analytic  family of $\phi$-invariant metrics representing 
the constructed conformal classes $[g_t]$, with the desirable property that  the 
 scalar curvature is of a fixed sign  for each  $t$.  

But what {\sl is} the  sign of the scalar curvature? To answer this, first 
 observe that, for each $t\in (0,1)$, 
the universal cover of $Y$ 
can naturally be realized as an open set of $S^4$, namely the union $\Omega_t$ of
all translates of $D_t$ via elements of the group  generated by the M\"obius transformations
$q\mapsto t^2 (q+i)^{-1} +i $ and $q\mapsto -t^2 (q -\sqrt{3})^{-1}-\sqrt{3}$
of $S^4 = \HH\PP_1$. In other words, 
$$(S^1 \times S^3)\# (S^1 \times S^3)= \Omega_t/ \ZZ \ast \ZZ,$$
where $\Omega_t$ is the region of discontinuity of the Kleinian group
$$\ZZ\ast \ZZ \subset PSL(2, \CC) \subset PGL (2 , \HH )$$ generated by 
$$
\frac{1}{t}
\left[ 
\begin{array}{rc}
1&i-it^2\\
-i&1
\end{array}
\right]
~~\mbox{and}~~
\frac{1}{t}
\left[ \begin{array}{cc}
\sqrt{3}&t^2-3\\
-1&\sqrt{3}
\end{array}\right] ~~
\in ~ SL(2, \CC). 
$$
These Kleinian groups are of the special type known as {\em Schottky groups}  \cite{maskit}.
Henceforth,  $D_t$ is to   be understood as a fundamental domain for the 
corresponding group action.

The complement $\Lambda_t= S^4 -\Omega_t$ 
of the region of discontinuity is 
called the {\em limit set} of the  group action.
If we think of $S^4$ as the boundary of the $5$-disk, on
whose interior $PGL (2, \HH ) = SO^\uparrow (5,1)$ acts
by isometries of the hyperbolic metric, then the limit set may also be 
characterized as the accumulation points of the orbit of any point in the open $5$-ball. 
 Since we have arranged for each of our
subgroups of $PGL (2, \HH ) = SO^\uparrow (5,1)$ to actually lie in 
$PSL (2, \CC ) = SO^\uparrow (3,1)$, it follows that we  have $\Lambda_t \subset 
\CC\PP_1 \subset \HH \PP_1$. This  will later allow us  to   
use planar diagrams to understand the structure of the limit set.
 
For our purposes, 
the ultimate utility of the Kleinian point of view  stems from a
 remarkable result of Schoen and Yau \cite{schyaudim}  that 
 relates  the scalar curvature of a uniformized 
 conformally flat manifold to the size  of
 the corresponding limit set. The form of this result we will use  is actually a  slight refinement 
 due to  Nayatani \cite{shinpat}:

\begin{lem}[Schoen-Yau, Nayatani] 
\label{schauen}
Let (M,[g]) be a compact, locally conformally flat 
$n$-manifold, $n\geq 3$,  which can be uniformized as $$M=\Omega/G,$$ where 
$G\subset SO^\uparrow (n+1, 1)$ is a Kleinian group
and where $\Omega\subset S^n$ is the region of discontinuity of $G$. Let 
$g\in [ g]$ be a  metric on $M$ in the fixed conformal class for which 
the scalar curvature $s$ does not change sign. Assume that  the limit set $\Lambda$ of $G$
is infinite, and let $\dim (\Lambda ) > 0$ denotes its 
Hausdorff dimension.
Then
\begin{eqnarray*}
s > 0 & \Longleftrightarrow & \dim (\Lambda ) < \frac{n}{2}-1 \\
s = 0 & \Longleftrightarrow & \dim (\Lambda ) = \frac{n}{2}-1 \\
s < 0 & \Longleftrightarrow & \dim (\Lambda ) > \frac{n}{2}-1  . 
\end{eqnarray*} 
\end{lem}

The original argument given by  Schoen and Yau is rather indirect, but 
Nayatani's proof actually constructs a particular  metric
for which the scalar curvature does not change sign; his conformal factor
is obtained by convolving an appropriate power of the Euclidean distance 
with  the Patterson-Sullivan measure of the limit set.  Because this construction is so natural
and canonical, 
it might seem 
 tempting to simply use Nayatani's algorithm to define our family of metrics $g_t$.
  We have  avoided doing so here, however, in order avoid the technical 
problem of proving that these metrics  depend analytically on the parameter $t$.

To prove Proposition \ref{cle}, we  now proceed by   
showing that $\dim (\Lambda_t) < 1$
for $t$ close to $0$, and that $\dim (\Lambda_t) > 1$ for $t$ close to $1$. 

Since the region of discontinuity $\Omega_t$ is the union  of all translates of $D_t$,
the limit set $\Lambda_t$ may be thought of as the intersection of a nested
sequence of balls in $\RR^4 = \HH$, where each of our original four balls
contains the reflections of the other three, each of these in turn contains another
three, and so forth. However, we have also observed that 
$\Lambda_t=\Lambda_t\cap \CC$, so the limit set may instead be 
thought of as the generalized Cantor set in $\CC$ given by 
the 
intersection of a  nested sequence of disks,
where,  in passing from one level to the next,  each disk is replaced 
by three smaller ones.
 Now    the  two generators of our Schottky group both  have derivatives satisfying  
$$ \left|\frac{d}{dz}\left(\pm\frac{t^2}{z-c}\pm c\right)\right| 
\leq  t^2 ~~~\mbox{whenever}~~ |z-c|\geq 1.$$
Since the  disks at the $k^{\rm th}$ level of the nesting are obtained by applying 
compositions of 
$k$ generators to one of the original $4$ disks, this implies that the  
 the disks at the $k^{\rm th}$ level  have Euclidean radius $< t^{2k}$. 
There are $4\cdot 3^k$ of these, 
  \begin{center}
\mbox
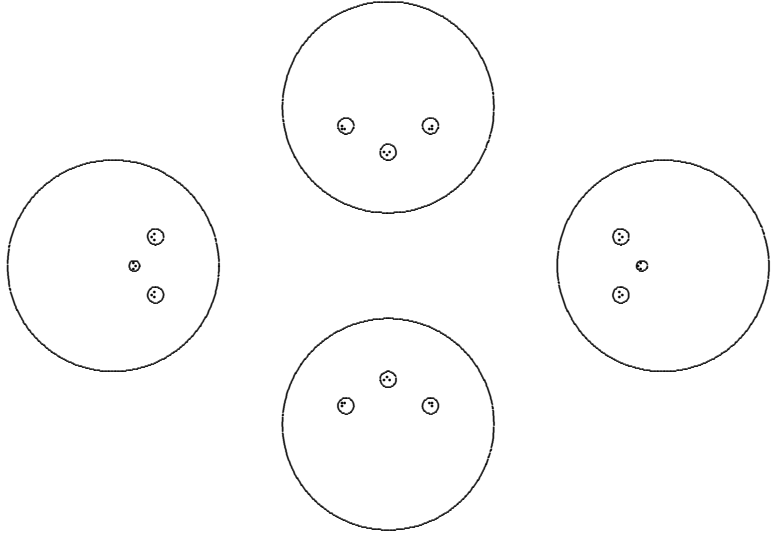
{
\beginpicture
\setplotarea x from -150 to 150, y from -100 to 100
 \circulararc 360 degrees from 64 0 center at  104 0
  \circulararc 360 degrees from -64 0 center at  -104 0
\circulararc 360 degrees from 0 20 center at  0 60
\circulararc 360 degrees from 0 -20 center at  0 -60
 \circulararc 360 degrees from 0 40 center at  0 43
  \circulararc 360 degrees from 0 -40 center at  0 -43
  \circulararc 360 degrees from 16 50 center at  16 53
   \circulararc 360 degrees from 16 -50 center at  16 -53 
   \circulararc 360 degrees from -16 50 center at  -16 53 
   \circulararc 360 degrees from -16 -50 center at  -16 -53
    \circulararc 360 degrees from 88 8 center at  88 11 
    \circulararc 360 degrees from -88 8 center at  -88 11 
    \circulararc 360 degrees from -88 -8 center at  -88 -11 
    \circulararc 360 degrees from 88 -8 center at  88 -11 
    \circulararc 360 degrees from 94 0 center at  96 0 
       \circulararc 360 degrees from -94 0 center at  -96 0 
       \put {\circle*{1}} [B1] at 1 43
          \put {\circle*{1}} [B1] at -1 43
             \put {\circle*{1}} [B1] at 0 42
               \put {\circle*{1}} [B1] at 1 -43
                \put {\circle*{1}} [B1] at -1 -43
             \put {\circle*{1}} [B1] at 0 -42
          \put {\circle*{1}} [B1] at 17 52
             \put {\circle*{1}} [B1] at 16 52
              \put {\circle*{1}} [B1] at 17 53
              \put {\circle*{1}} [B1] at 17 -52
             \put {\circle*{1}} [B1] at 16 -52
              \put {\circle*{1}} [B1] at 17 -53
              \put {\circle*{1}} [B1] at -17 52
             \put {\circle*{1}} [B1] at -16 52
              \put {\circle*{1}} [B1] at -17 53
              \put {\circle*{1}} [B1] at -17 -52
             \put {\circle*{1}} [B1] at -16 -52
              \put {\circle*{1}} [B1] at -17 -53
                \put {\circle*{1}} [B1] at 88 10
             \put {\circle*{1}} [B1] at 88 12
              \put {\circle*{1}} [B1] at 89 11
               \put {\circle*{1}} [B1] at 88 -10
             \put {\circle*{1}} [B1] at 88 -12
              \put {\circle*{1}} [B1] at 89 -11
 \put {\circle*{1}} [B1] at -88 10
             \put {\circle*{1}} [B1] at -88 12
              \put {\circle*{1}} [B1] at -89 11
 \put {\circle*{1}} [B1] at -88 -10
             \put {\circle*{1}} [B1] at -88 -12
              \put {\circle*{1}} [B1] at -89 -11
\put {\circle*{1}} [B1] at 95 0
             \put {\circle*{1}} [B1] at 96 1
              \put {\circle*{1}} [B1] at 96 -1
              \put {\circle*{1}} [B1] at -95 0
             \put {\circle*{1}} [B1] at -96 1
              \put {\circle*{1}} [B1] at -96 -1
\endpicture
}
\end{center}
so the $d$-dimensional Hausdorff measure of $\Lambda_t$  
 is  less than a constant times $(3t^{2d})^k$ for all $k$, and vanishes if
 $\log 3 + 2d \log t   < 0$. 
It therefore follows that 
$$\dim (\Lambda_t ) \leq  - \frac{\log 3}{ 2\log t } .$$
In particular,  $\dim (\Lambda_t ) < 1$ for all $t \in (0,1/2]$, and for these values of  $t$
our $\phi$-invariant conformally flat metrics $g_t$
will have $s > 0$.  The interested reader may enjoy the exercise of constructing
explicit  choices of  $g_t$ with $s > 0$  when $t$ is extremely small, and
comparing the results obtainable in this way with the predictions of the above limit-set argument. 

Next, we need to show that $\dim (\Lambda_t ) > 1$ when $t$ is sufficiently
close to $1$. To see this, first consider the Kleinian group $G\subset PSL(2, \CC)$
generated by 
$$
\left[ 
\begin{array}{rc}
1&0\\
-i&1
\end{array}
\right]
~~\mbox{and}~~
\left[ \begin{array}{cc}
\sqrt{3}&-2\\
-1&\sqrt{3}
\end{array}\right] ~~
\in ~ SL(2, \CC). 
$$
which is  the limiting case of our  construction that arises by na\"{\i}vely setting  $t=1$. 
We can still construct a fundamental domain $D$ domain for this action 
as the complement of four balls, but  certain pairs our four  balls in $\HH$ will 
now have a boundary point
in common.
The corresponding system of nested disks in  $\CC$ then contains 
two `bracelets' of disks arranged around a pair of circles in the plane: 
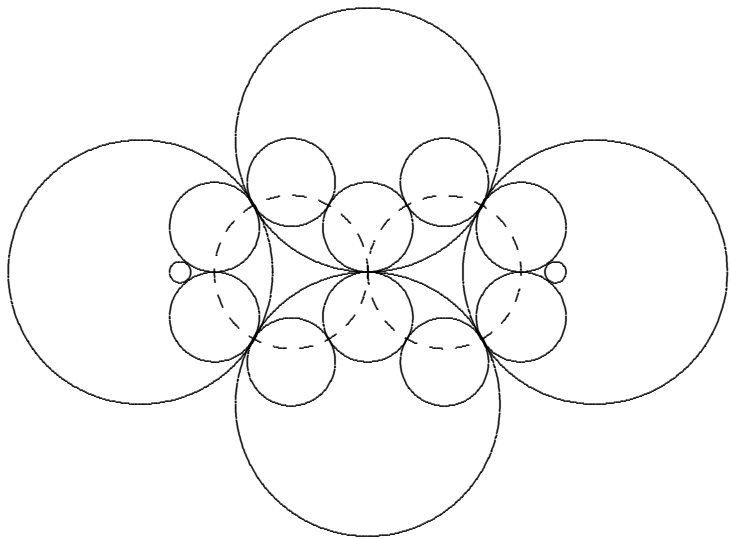
\begin{center}
\mbox{
\beginpicture
\setplotarea x from -150 to 150, y from -105 to 100
 \circulararc 360 degrees from 36 0 center at  86 0
  \circulararc 360 degrees from -36 0 center at  -86 0
\circulararc 360 degrees from 0 0 center at  0 50
\circulararc 360 degrees from 0 0 center at  0 -50
\circulararc 360 degrees from 0 0 center at  0 17
\circulararc 360 degrees from 0 0 center at  0 -17
\circulararc 360 degrees from 43 25 center at  29 34
\circulararc 360 degrees from 43 -25 center at  29 -34
\circulararc 360 degrees from -43 25 center at  -29 34
\circulararc 360 degrees from -43 -25 center at  -29 -34
\circulararc 360 degrees from 43 25 center at  58 17
\circulararc 360 degrees from -43 25 center at  -58 17
\circulararc 360 degrees from 43 -25 center at  58 -17
\circulararc 360 degrees from -43 -25 center at  -58 -17
\circulararc 360 degrees from 67 0 center at  71 0
\circulararc 360 degrees from -67 0 center at  -71 0
{\setdashes
\circulararc 360 degrees from 0 0  center at  29 0
\circulararc 360 degrees from 0 0  center at -29 0
}
\endpicture
}
\end{center}
As one passes from one level of the nesting to the next, both of these circles continue 
to be completely covered by bracelets of smaller and smaller disks. 
Hence these two circles are both  contained in the  
limit set of $G$. But, by a result of Bishop and Jones  \cite{bishopdim},  this implies 
 that the Hausdorff dimension of the corresponding limit set  
 must be strictly greater than one: 
 
 \begin{lem}[Bishop-Jones]
 Let $G\subset PSL (2, \CC )$ be a  finitely generated
 Kleinian group with infinite limit set $\Lambda\subset {\mathbb C \mathbb P}_1$. 
 If $\Lambda$ is not totally disconnected, then either $\Lambda$ is  
a single geometric circle, or else $\dim (\Lambda ) > 1$. 
 \end{lem}
 
 Now, to clinch the argument, we would like to somehow use this  to estimate the Hausdorff dimension
of  $\Lambda_t$
as $t\to 1$. 
 A second general result of Bishop and Jones \cite{bishopdim} provides 
 the machinery  needed to do  this: 
 
  \begin{lem}[Bishop-Jones]
 Let $G\subset PSL (2, \CC )$ be a  finitely generated
 Kleinian group, and suppose that $G$ is  the limit of subgroups
 $G_j\subset PSL (2, \CC )$, in the sense of convergence of  a set of generators. 
 If $\Lambda$ is the limit set of $G$, and if $\Lambda_j$ is the 
 limit sets of   $G_j$, 
then $\liminf \dim (\Lambda_j)  \geq \dim (\Lambda )$. 
 \end{lem}
 In our case, it follows that there is an $\epsilon > 0$ such that 
 $$\dim (\Lambda_t ) > 1~~~\mbox{whenever}~~~ t > 1- \epsilon,$$
 since otherwise there would exist a sequence $t_j\nearrow 1$
 with $\dim (\Lambda_{t_j}) \leq 1$, and hence with 
 $\liminf \dim (\Lambda_{t_j}) \leq 1 < \dim (\Lambda )$, in contradiction to the lemma.  
For $t\in (1-\epsilon , 1 )$, the corresponding $\phi$-invariant 
conformally flat metrics $g_t$ therefore 
have  $s< 0$.

To wrap up our proof of Proposition \ref{cle}, it just remains to  show that 
\begin{equation}
\label{eigenzero}
B=\Big\{ t \in (0,1) ~\Big|~ g_t \mbox{ has } \ker (\Delta + \frac{s}{6})\neq  0\Big\},
\end{equation}
is a discrete set, so that the subset $A\subset B\cap  [\frac{1}{2}, 1-\epsilon]$
defined by 
\begin{equation}
\label{scalzero}
A=\Big\{ t \in (0,1) ~\Big|~ g_t \mbox{ has } s\equiv 0\Big\}
\end{equation}
is consequently finite. In principle this could  again be done by appealing to
the perturbation theory of the
spectrum of the Yamabe Laplacian, but, just for fun, let us give a twistorial proof, in the
spirit of \S \ref{greatsmall}. Indeed, 
 the Penrose transform tells us  that we can re-express (\ref{eigenzero}) as  
$$B=\Big\{ t \in (0,1)~\Big|~ H^1(Z_t , {\mathcal O} (K^{1/2}))\neq 0\Big\}, $$
where $Z_t$ is the twistor space of $(Y,g_t)$.  However, $Z_t$
is constructed by taking an open set in  $\CC\PP_3$ (namely the the inverse image
of some open neighborhood of $D_t\subset \HH \PP_1$ 
via the twistor projection) and making identifications
using the two biholomorphisms given by the  $PSL(4, \CC)$-transformations 
arising from the generators of our Schottky group
$\ZZ \ast \ZZ$ via the inclusions 
$$PSL (2, \CC ) \hookrightarrow PGL (2, \HH ) \hookrightarrow PSL (4, \CC ).$$
But our generators depend algebraically on $t$, so we may  extend our 
construction of the $Z_t$ to $t$ in an open neighborhood $\mathcal U$ of $(0,1) \subset \CC$,
giving us an analytic family of complex $3$-folds. The semi-continuity principle \cite{bast}
then 
implies that the set of $t\in {\mathcal U}$ for which $H^1 (Z_t , {\mathcal O} (K^{1/2}))\neq 0$
is closed in the analytic Zariski topology; in other words, it is either discrete, or else 
is all of $\mathcal U$. Since $\Delta + \frac{s}{6}$ is a positive operator for 
$t$ small, this shows that the set $B$ defined by (\ref{eigenzero}) is 
discrete. The compact set $A\subset B$ defined by (\ref{scalzero}) is thus finite, as claimed. 

It follows that the element of $A$ defined by 
$$t_0= \sup \{ t\in (0,1) ~|~ g_t \mbox { has } s> 0\}$$
has a neighborhood $(t_0 - \varepsilon , t_0+ \varepsilon )\subset (0,1)$ which does not meet 
$A - \{ t_0\}$. For $t$ in this neighborhood, the scalar curvature $s$ of 
$g_t$ then has the same sign as $t_0-t$, and  our proof of 
Proposition \ref{cle} is therefore done.

\section{A Vanishing Theorem}
\label{vanish}

At this point, we have constructed an interesting family of locally conformally flat
metrics on the orbifold $X=Y/\ZZ_2$, where $Y=2(S^1 \times  S^3 )$.
However,  
our aim is to eventually  smooth the orbifold
singularities of  these metrics in order to 
produce a similar family 
of anti-self-dual metrics on a simply connected manifold.
To carry this out, we will need to know that the Kodaira-Spencer deformation theory
is unobstructed for the  corresponding family of twistor spaces. 
In fact, 
the relevant vanishing theorem  easily follows from a
decade-old 
unpublished paper of Eastwood and Singer \cite{eastsing},
whose  beautiful  ideas will be given a
 self-contained exposition  in this section.

Let $(Y,g)$ be an oriented,  locally conformally flat  Riemannian $4$-manifold,
and let $\mathscr C$ denote the complete presheaf of {\em conformal Killing fields}
on $(Y,g)$, defined by setting  
$${\mathscr C}_U= \{ v\in {\mathcal E}_U(TY)~|~\pounds_vg\propto g\}$$
for any open set $U\subset Y$; here, as throughout, 
${\mathcal E}$ is used to indicate the  $C^{\infty}$ sections of a given
vector bundle. Now observe that there is a rank-$15$ vector bundle
$F\to Y$, equipped with a flat connection $\triangledown$, such 
that $\mathscr C$  is the sheaf of parallel sections of $(F, \triangledown )$. 
Indeed, if $U\subset Y$ is any {\em simply connected} open set, 
then we can conformally immerse $U$ onto an open subset of $S^4$
by means of the {\em developing map} \cite{kuiper}, and 
${\mathscr C}_U$ is thereby identified 
with the $15$-dimensional space ${\mathfrak s\mathfrak o}(5,1)$ of global
conformal Killing fields on the round $4$-sphere $S^4$. As we pass from one
such choice of $U$ to another, these identifications will be related to one another 
by  elements of $SO^\uparrow(5,1)$, acting on  ${\mathfrak s\mathfrak o}(5,1)$
via the adjoint representation. These elements of  $SO^\uparrow(5,1)$ are exactly the
transition of functions of $F$, relative to  a collection of local trivializations of 
$F$ in which the flat connection $\triangledown$ has vanishing connection $1$-forms. 

Now we could certainly construct a fine resolution of $\mathscr C$ by just considering
the $F$-valued differential forms on $Y$,  but this would involve  
using vector bundles   of rather high rank. 
A  more efficient  resolution was first discovered by 
Gasqui and Goldschmidt \cite{gasgold} using Spencer cohomology, and later
rediscovered  by Eastwood
and Rice \cite{rice} in the setting of Bernstein-Gelfand-Gelfand resolutions. This 
resolution 
takes the form
$$0\to{\mathscr C}\to 
{\mathcal E}(TY)\stackrel{L_0}{\to}{\mathcal E}(\odot^2_0\Lambda^1)
\stackrel{L_1}{\to}{\mathcal E}\Big(
\begin{array}{c}\odot^2_0\Lambda^+
\\
\oplus
\\
\odot^2_0\Lambda^-
\end{array}
\Big)
\stackrel{L_2}{\to}{\mathcal E}(\odot^2_0\Lambda^1)\stackrel{L_3}{\to}
{\mathcal E}(TY) \to 0
$$
where $\odot^2_0$ indicates the trace-free symmetric-tensor-product of a vector bundle
with itself. Here $L_0$ is the first-order operator 
$$L_0(v)= \mbox{trace-free part of } {\pounds}_vg $$ 
which measures the way the conformal class $[g]$  is distorted by the 
the flow of a given vector field. 
The next step in the sequence is the linearization  $L_1=DW$
of the Weyl curvature tensor; and for our purposes, it will be important to
recognize this second-order differential operator can be decomposed
as $DW= DW_++ DW_-$, where the operators
\begin{eqnarray*}
DW_+:  & {\mathcal E}(\odot^2_0\Lambda^1)\to  & {\mathcal E}(\odot^2_0\Lambda^+ )\\
DW_-:  & {\mathcal E}(\odot^2_0\Lambda^1)\to  & {\mathcal E}(\odot^2_0\Lambda^- )
\end{eqnarray*}
are the linearizations 
$$DW_\pm (h) = \frac{d}{dt}W_\pm(g+th)\Big|_{t=0}$$
 of the self-dual and anti-self-dual Weyl curvatures. 
 The next step is again a 
 second-order operator, and is given by  
$$L_2= (DW_+)^*-(DW_-)^*.$$
The sequence  then culminates with   the first-order operator 
$$L_3=L_0^*.$$
Clearly, all of these operators are conformally invariant, provided that each 
bundle in the complex is given the   correct conformal weight.

Since each of the sheaves in the Gasqui-Goldschmidt resolution is fine, and hence acyclic,
the abstract de Rham theorem \cite{wells} 
immediately tells us that the sheaf
 cohomology of  $Y$ with coefficients in $\mathscr C$ is exactly the  cohomology of 
 the corresponding 
complex of global sections:
 $$H^p(Y, {\mathscr C}) = \frac{\ker L_{p}}{\mbox{im }  L_{p-1}}.$$
 However,  the Gasqui-Goldschmidt resolution is also  
 an elliptic complex; thus,   provided $Y$ is {\em compact},
 we have   \cite{gasgold2}
   $$H^p(Y, {\mathscr C}) =\ker L_{p}\cap \ker  L_{p-1}^*,$$ 
   by a generalized form of the Hodge theorem. 
   Since 
$$DW_\pm = \frac{1}{2}(L_1^* \pm L_2),$$
this immediately gives us  the following result:

\begin{prop} \label{hodge}
Let $(Y,g)$ be  any compact, oriented, locally conformally flat 
$4$-manifold. Then 
\begin{eqnarray*}
H^2(Y,{\mathscr C})  & = & \ker  (DW_+)^* \oplus \ker  (DW_-)^* \\
 & \subset & {\mathcal E}_Y(\odot^2_0\Lambda^+)~\oplus  ~{\mathcal E}_Y(\odot^2_0\Lambda^+).
\end{eqnarray*}
\end{prop}

Using this key observation, it is now  easy to deduce the desired  vanishing result:
 
\begin{thm}[Eastwood-Singer] \label{singeast}
Let $g$ be {\em any} conformally flat metric on the oriented $4$-manifold 
$$Y=k(S^1\times S^3)=\underbrace{(S^1\times S^3)\# \cdots \# (S^1\times S^3)}_k,$$
$k\geq 1$, 
and let $Z$ be the   twistor space of $(Y,g)$. 
Then 
$$DW_+: {\mathcal E}(\odot^2_0\Lambda^1) \to  {\mathcal E}(\odot^2_0\Lambda^+ )
 $$
is surjective on $(Y,g)$, and 
$$H^2(Z, {\mathcal O}(TZ))=0.$$
\end{thm}
\begin{proof}
By Serre duality, $H^2(Z, {\mathcal O}(TZ))$  is the dual 
of  $H^1(Z, \Omega^1(K))$.
However, the latter sheaf cohomology group corresponds, via the Penrose transform 
\cite{bailsing}, 
to $\CC\otimes \ker (DW_+)^*$. By Proposition \ref{hodge}, it therefore suffices to 
show  that $H^2(Y, {\mathscr C}) = 0$, where ${\mathscr C}$ is once again the sheaf of
local conformal Killing fields of $[g]$. 

Now $Y$ can be obtained from $S^4$ 
by replacing $k$ pairs of  balls with $k$ tubes modeled on $S^3 \times \RR$. 
This allows us to  
express $Y$ as the union 
$$Y = U \cup V$$
of  open sets 
$$U = S^4 -\{ p_1, \ldots, p_{2k}\}$$
and 
$$V \approx \underbrace{(S^3 \times \RR)\sqcup   \cdots \sqcup    (S^3 \times \RR)}_k $$
such that 
$$U\cap V \approx \underbrace{(S^3 \times \RR)\sqcup   \cdots \sqcup    (S^3 \times \RR)}_{2k}.$$
We may thus proceed by examining  the Mayer-Vietoris sequence 
\begin{equation}
\label{mv}
 \to H^1(U\cap V,  {\mathscr C})
\to H^2(U\cup V, {\mathscr C}) \to  H^2(U,  {\mathscr C})\oplus
H^2(V, {\mathscr C})\to \cdots 
\end{equation}
Indeed, 
notice that  $V$ and $U\cap V$ are homotopy equivalent to  disjoint unions of $3$-spheres,
while $U$ is homotopy equivalent to a bouquet of $3$-spheres. In particular,  each of these
sets is  a disjoint union of  simply connected spaces.  Since $\mathscr C$
is the sheaf of covariantly constant  sections of a flat rank-$15$
vector bundle $(F,\triangledown )$, 
the restriction of $\mathscr C$ to any of these open sets may  be identified 
with the constant sheaf $\RR^{15}$, and  the relevant sheaf cohomology therefore amounts to   singular cohomology with coefficients in the Abelian group $\RR^{15}$. 
By the homotopy invariance of  singular cohomology, we thus have 
\begin{eqnarray*}
H^1(U\cap V,  {\mathscr C})&\cong & H^1(\underbrace{S^3\sqcup   \cdots \sqcup   S^3}_{2k}, \RR^{15}) =0 \\
H^2(U,  {\mathscr C})& \cong & H^2(\underbrace{S^3\vee
 \cdots \vee S^3}_{2k-1}, \RR^{15}) =0  \\
H^2( V,  {\mathscr C})& \cong & H^2(\underbrace{S^3\sqcup   \cdots \sqcup   S^3}_{k}, \RR^{15}) =0, 
\end{eqnarray*}
and (\ref{mv})  therefore tells us that  
$$H^2(Y, {\mathscr C})=H^2(U\cup V, 
{\mathscr C})=0,$$  
as claimed. 
\end{proof}

It is perhaps worth remarking that, for {\em any} 
compact oriented locally conformally flat $4$-manifold $(Y,g)$, 
one may use  the index theorem to show that 
 $\ker (DW_+)^*$ and 
$\ker (DW_-)^*$   have  the same dimension. 
Proposition \ref{hodge} and  the Penrose transform
therefore imply that 
$$\dim_\RR H^2(Y, {\mathscr C}) = 2 \dim_\CC H^2 (Z, {\mathcal O}(TZ)).$$
Thus the vanishing of $H^2(Y, {\mathscr C})$ is actually
  necessary, as well as  sufficient, 
for the deformation theory of $Z$ to be unobstructed.

\section{Existence Results}

We will now assemble the results of the last several sections  into 
a proof of Theorem \ref{able}.

\begin{prop}\label{change}
For any integer $k\geq 6$, the  connected sum 
$$k\overline{\bcp}_2=\underbrace{\overline{\bcp}_2\# \cdots \# \overline{\bcp}_2}_k$$ 
admits
a real-analytic $1$-parameter family of anti-self-dual  conformal metrics $[g_t]$, 
$t\in [a,b]$,
such that $[g_a]$ contains a metric with $s> 0$ everywhere, while
$[g_b]$ contains a metric with $s<0$ everywhere.
\end{prop}
\begin{proof}
We once again let $Y=(S^1\times S^3) \# (S^1\times S^3)$, and let 
$X=Y/\ZZ_2$.  Equip $X$ with a conformally-flat scalar-flat  orbifold metric $g_0$
which belongs to a real-analytic $1$-parameter family of conformally flat metrics $g_t$, 
$t\in (-\varepsilon , \varepsilon )$
such that  each $g_t$  has $s>0$ when $t<0$,  $s<0$ when 
$t>0$, and 
$\ker (\Delta + s/6) =0$ for all $t\neq 0$; for example, 
Proposition \ref{cle} constructs such a family,
after replacing $t$ with $t-t_0$ for some $t_0$. 
 Since   
Theorem \ref{singeast} guarantees that the twistor space $Z$ of $(Y, g_0)$ 
has $H^2(Z, {\mathcal O}(TZ))=0$, 
Theorem \ref{thumbo} therefore tells us that, for any integer 
$\ell \geq 0$,
 there exists a connected $2$-parameter family of anti-self-dual conformal metrics 
 $[g_{ ( {\mathfrak u}_{1} , {\mathfrak u}_{2} ) }]$, $({\mathfrak u}_1,{\mathfrak u}_2)\in 
 (0 , \epsilon ) \times (-\epsilon , \epsilon )$ 
on $\tilde{X}\# \ell \overline{\bcp}_2$ which arises from a $2$-parameter standard deformation of the 
singular space $Z_0=[\tilde{Z}/\ZZ_2]\cup \tilde{Z}_{EH}\cup \tilde{Z}_{FS}$;
here 
 $\tilde{X}$ is once again 
 the oriented $4$-manifold obtained from the orbifold $X$ by replacing 
each singular point of $X$ with a $2$-sphere of self-intersection $-2$.
Moreover this $2$-parameter standard deformation can be 
chosen so that its restriction to ${\mathfrak u}_2=t$ is a $1$-parameter 
standard deformation of the complex  $3$-fold with normal crossings 
arising from 
$(X, [ g_t ])$,  for every  real number $t$ in a neighborhood of  $0$. 
Let $a< 0$ and $b>0$ be choices of $t$ in this neighborhood. Theorems \ref{neg}
and \ref{pos} then tell us that for any sufficiently small $w >0$, the conformal class 
$[g_{(w, a)}]$ contains a metric with $s > 0$, while the conformal class $[g_{(w, b)}]$
contains a metric with $s < 0$. Thus $[g_t]:=[g_{(w,t)}]$, $a\leq t\leq b$,  is a
family of ASD conformal metrics on $\tilde{X}\# \ell \overline{\bcp}_2$ with the
desired scalar-curvature behavior .

 It remains only to unmask the identity of the manifold $\tilde{X}\# \ell \overline{\bcp}_2$.
 To do this, first notice we may cut up $Y$ into three punctured $4$-spheres
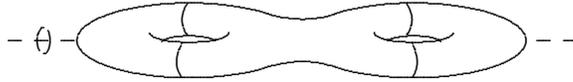
\begin{center}
\mbox{
\beginpicture
\setplotarea x from 40 to 290, y from -5 to 65
\ellipticalarc axes ratio 3:1  270 degrees from 150 40
center at 120 30
\ellipticalarc axes ratio 3:1  -270 degrees from 175 40
center at 205 30
\ellipticalarc axes ratio 4:1 -180 degrees from 135 33
center at 120 33
\ellipticalarc axes ratio 4:1 145 degrees from 130 30
center at 120 29
\ellipticalarc axes ratio 4:1 180 degrees from 190 33
center at 205 33
\ellipticalarc axes ratio 4:1 -145 degrees from 195 30
center at 205 29
\ellipticalarc axes ratio 1:2 140 degrees from 65 35
center at 65 30
\ellipticalarc axes ratio 1:2 140 degrees from 65 25
center at 65 30
{\setquadratic 
\plot 150 40    163 38    175 40   /
\plot 150 20    163 22    175 20  /
}
{\setquadratic 
\plot 116 32   114 38  116 44 /
\plot 199 32   201 38   199 44 /
\plot 114 29   112 22    114 16 /
\plot 201 29   203 22    201 16 /
}
{\setlinear
\setdashes 
\plot  45 30 70 30   /
\plot  245 30 265 30   /
}
\endpicture
}
\end{center}
in a manner 
which is compatible with the involution. Thus 
  $X$ can be expressed as a connected sum 
$$X = V\# V \#V$$
of three copies of the orbifold $S^4/\ZZ_2$, 
where the connect sum is  
carried out in the vicinity of  {\em non-singular} points of $V$, 
and where the $\ZZ_2$ acts on $S^4 \subset \RR^5$ by reflection through an axis. Hence 
$$\tilde{X} = \tilde{V}\# \tilde{V} \# \tilde{V} , $$
where $\tilde{V}$ is obtained from $V$ by replacing the two isolated singularities of 
$V=S^4/\ZZ_2$ by $2$-spheres of self-intersection $-2$. However,
$\overline{\bcp}_2 \# \overline{\bcp}_2$  is diffeomorphic to 
$\tilde{V}$. Indeed, 
if $E_1$ and $E_2$ are the standard $2$-spheres of self-intersection $-1$ 
in $\overline{\bcp}_2 \# \overline{\bcp}_2$, then $E_1+E_2$ and $E_1-E_2$
are  represented by disjoint embedded $2$-spheres of self-intersection $-2$, and 
$\overline{\bcp}_2 \# \overline{\bcp}_2$ is obtained by gluing tubular neighborhoods 
of these two $2$-spheres along their boundaries. 
We therefore have $X\approx  6\overline{\bcp}_2$, and hence 
$M\approx k\overline{\bcp}_2$, where $k=6+\ell$. 
\end{proof}

As a corollary, we now obtain one of the  central results of this paper:

\begin{thm} \label{capstone}
 For any integer $k\geq 6$, the  connected sum $k\overline{\bcp}_2$ admits
scalar-flat anti-self-dual  metrics. 
\end{thm}
\begin{proof}
Consider the smooth family of conformal classes $[g_t]$ constructed in 
Proposition \ref{change}, and let $h_t\in [g_t]$ be any smooth family of
metrics representing these conformal classes. Let $\lambda_t$
denote the smallest eigenvalue of the Yamabe Laplacian $(\Delta + s/6)$
for the metric $h_t$. Then $\lambda_t$ is a continuous function of $t$.
But Proposition \ref{change} tells us that $\lambda_a > 0$, whereas
$\lambda_b < 0$. By continuity, there is thus some $c\in [a,b]$ for which
$\lambda_c=0$. Let $u$ be a  unit-integral eigenfunction of the 
Yamabe Laplacian $(\Delta + s/6)$ of $h_c$ with eigenvalue $\lambda_c=0$. 
By the minimum principle, $u$ is a positive function. Thus $g=u^2h_c$
is  a scalar-flat anti-self-dual metric on $k\overline{\bcp}_2$
for the given value of $k\geq 6$. 
\end{proof}

Now this by no means represents  the first construction ever  of SFASD metrics on  simply connected 
compact $4$-manifolds.  However, all the previous results depended on
an essentially different idea: namely, that 
 any K\"ahler metric on a complex surface with $s\equiv 0$
is SFASD.  Through this observation, 
Yau's  existence theorem for  Ricci-flat K\"ahler metrics on 
$K3$ surfaces  \cite{yauma} provided a crucial early  family of examples which largely drove
the subsequent development of the entire subject. Much later, 
the present author and his collaborators showed \cite{klp}  that 
$\bcp_2\#k \overline{\bcp}_2$ admits scalar-flat anti-self-dual metrics if
$k\geq 14$. The proof of this last result depends on a  refinement of Theorem \ref{thmb}, 
set up so that the constructed twistor spaces carries a  special divisor
whose existence implies that the ASD conformal class contains a K\"ahler metric. 
For a related  re-proof of the existence of Calabi-Yau metrics on $K3$, see \cite{lebsing2}.

Putting these previous results  together with Corollary \ref{capstone}, 
we have thus proved Theorem \ref{able}:

\begin{thm}
The following   smooth $4$-manifolds   admit scalar-flat anti-self-dual   
 metrics:
 \begin{romenumerate}
 \item   $k \overline{\bcp}_2$, for  every $k  \geq 6$; 
\item  $\bcp_2\#k \overline{\bcp}_2$, for every  $k \geq 14$; and  
\item  $K3$.
\end{romenumerate}
In particular, each of these simply connected compact 
$4$-manifolds admits optimal metrics; and any optimal metric
on any one of them is SFASD. 
\end{thm}

\section{Non-Existence Results}
\label{nonex}

We have now seen that many simply connected $4$-manifolds admit non-Einstein 
optimal metrics. However, related ideas will now allow us to   show 
there are also many simply connected $4$-manifolds
which do {\em not} admit optimal metrics. 
To see this, we begin by introducing a new concept:

\begin{defn}
Let $M$ be a smooth compact oriented $4$-dimensional manifold (respectively, orbifold). 
We will say that $M$ admits an {\em anorexic sequence} if there is a sequence
$g_j$ of smooth Riemannian metrics (respectively, orbifold metrics) on $M$
for which $\int s^2d\mu \to 0$ and $\int |W_+|^2 d\mu \to 0$. 
\end{defn} 

When a manifold admits such a sequence, we then know the value of ${\mathcal I}_{\mathcal R}(M)$,
and stand a very good chance of determining whether it admits an optimal metric:

\begin{lem}\label{skinny}
Let $M$ be a smooth compact oriented $4$-manifold which admits an 
anorexic sequence. Then any optimal metric on $M$ is SFASD. 
Moreover, 
$${\mathcal I}_{\mathcal R}(M)= -8\pi^2(\chi + 3\tau ) (M).$$
\end{lem}
\begin{proof} Recall that  equation (\ref{knut}) tells us that 
$${\mathcal K}(g) = -8\pi^2 (\chi + 3\tau )(M) + 2\int_M \left(\frac{s^2}{24}+2|W_+|^2\right) d\mu_g ~. $$
If there is an anorexic sequence, the infimum of the right-hand side is thus obtained by 
dropping the curvature integral. Moreover, a metric minimizing $\mathcal K$ would necessarily
have $s\equiv 0$ and $W_+\equiv 0$.
\end{proof}

Now imagine a curvaceous  young $4$-manifold who, bedazzled by  
 the glamorous starlets with optimal metrics she has been reading about
 in the tabloids, suddenly decides to go on a starvation diet to get rid
 of  all that unwanted curvature. 
If she has the wrong body type, this misguided procedure will be dangerous to
her health, and she will  
merely succeed in putting herself in the 
hospital:

\begin{prop}\label{boney}
Let $M$ be a  smooth compact  oriented $4$-manifold which admits an 
anorexic sequence. Then $M$ does not admit an optimal metric if 
\begin{itemize}
\item $b_+(M) \geq 4$; or 
\item  $b_+(M) =2$; or 
\item $b_+(M) =3$, $\pi_1(M) =0$,  and $M$ is not diffeomorphic to $K3$; or 
\item $b_+(M) =1$, $\pi_1(M) =0$,  and $M$ is not diffeomorphic to  $\bcp_2 \# k \overline{\bcp}_2$
for some $k\geq 10$. 
\end{itemize}
\end{prop}
\begin{proof}
By Lemma \ref{skinny}, an optimal metric on such an $M$ would necessarily be SFASD.
However, 
Propositions \ref{form} and \ref{class} show that we would then obtain a contradiction 
in any of the scenarios considered above. \end{proof}

Now we come to our main method of 
construction \cite{lky}:

\begin{lem}\label{scrunch}
Let $X$ be an oriented compact $4$-dimensional orbifold with only isolated  singularities
modeled on $\RR^4/\ZZ_2$. Let $\tilde{X}$ be the smooth oriented $4$-manifold obtained 
by replacing each singular point by a $2$-sphere of self-intersection $-2$. If
$X$ admits an anorexic sequence,  then so does  $\tilde{X}$. Moreover, if there is an anorexic sequence
on $X$ with the property that $\int |r|^2d\mu \to 0$, then $\tilde{X}$ also admits an anorexic
sequence with this property. 
\end{lem}
\begin{proof}
If $h$ is an orbifold metric on $X$ with $\int s^2 d\mu < \varepsilon$ and 
$\int |W_+|^2 d\mu < \varepsilon$, we will  show that $\tilde{X}$ has a a metric
$g$ with $\int s^2 d\mu < 2\varepsilon$ and 
$\int |W_+|^2 d\mu < 2\varepsilon$. Moreover, 
 if  $h$ also has $\int |\mathring{r}|^2 d\mu  < \varepsilon$, then we will be able to arrange for 
 $g$ to also satisfy  $\int |\mathring{r}|^2 d\mu  < 2\varepsilon$.
 
 To do this, we choose geodesic normal coordinates around a given orbifold point,
 so that we have
 $$h=  \left[\sum_j (dx^j)^2\right]+\alpha $$
 in these coordinates, where $\alpha$ is  a  ($\ZZ_2$-invariant) smooth symmetric 
  tensor field on a neighborhood
 of the origin in $\RR^4$ with $|\alpha|< C|\vec{x}|^2$. 
 We would like to delete a ball of small radius $\rho$ around the origin, and
 glue in a copy of the Eguchi-Hanson metric on $T^*S^2$, with very small length scale. 
 Recall that the restriction of   the Eguchi -Hanson
 metric to  the complement of the zero section in $T^*S^2$ is isometric to the metric  
 $$g_{EH, \epsilon}= \frac{d\varrho^2}{1-(\frac{\epsilon}{\varrho})^4}+ \varrho^2\left( \sigma_1^2+\sigma_2^2+ 
\left[1-\left(\frac{\epsilon}{\varrho}\right)^4\right]\sigma_3^2\right)$$
 on $(\epsilon, \infty ) \times S^3/\ZZ_2$, where $\{ \sigma_j\}$ is 
the standard  left-invariant co-frame on $S^3/\ZZ_2 = SO(3)$;  the constant $\epsilon > 0$ is 
herein referred to as  the 
{\em length scale}.   Now, for any fixed $\rho >0$, this family of metrics  converges uniformly
in the $C^2$ topology to 
the Euclidean metric on the annulus $\varrho \in  [\rho /2 , \rho ]$ as $\epsilon \to 0$. 
If $\varphi : (0, \infty ) \to [0,1]$ is a bump function which is $\equiv 0$ on $(0,1/2]$
and $\equiv 1$ on $[1,\infty)$, then, for any fixed $\rho$, the metrics
$$g_{\epsilon, \rho} = \varphi  \left(
{
\textstyle \frac{\varrho}{\rho}
} \right)  h+ \left[1-\varphi  \left(
{\textstyle 
\frac{\varrho}{\rho} }
\right) \right] g_{EH, \epsilon}
$$
therefore converge  in the $C^2$ norm 
to 
$$g_{0, \rho}=\left[\sum_j (dx^j)^2\right]+ \varphi \left(
{\textstyle \frac{\varrho}{\rho}} \right) \alpha$$
 on the annulus $\varrho \in [\rho / 2 , \rho ]$,
and in particular the curvature tensors of these metrics converge uniformly in the annulus 
to the curvature of $g_{0, \rho}$.
On the other hand, since $\alpha$ is of magnitude $O(\varrho^2)$, the first and 
second coordinate partial derivatives
of $g_{0, \rho}$ are uniformly bounded as $\rho \to 0$. Thus we can 
choose a sequence of 
$(\epsilon_j , \rho_j) \to (0,0)$  such that the sectional curvatures of the 
metrics $g_j = g_{\epsilon_j , \rho_j}$ are uniformly bounded on the transition annuli
$\varrho \in [\rho_j/2 , \rho_j]$, while the volumes of these annuli simultaneously 
tend to zero. 
For $j$ far out in the sequence, the transition annulus therefore makes
a contribution to $\int s^2 d\mu$, $\int |\mathring{r}|^2d\mu$, or $\int |W_+|^2d\mu$
which is as small as we like --- for example, smaller than the given $\varepsilon$
divided by the number of orbifold singularities of  $X$. We now take  $g$
to be given by such a choice of $g_j$ in each glued region, $h$ on the complement
$X$ minus a collection of balls or radius $\rho_j$ about its orbifold singularities, and equal to the 
Eguchi-Hanson metric with length scale $\epsilon_j$   near the added $2$-spheres. 
Since the  Eguchi-Hanson metric
 has $r\equiv 0$ and $W_+\equiv 0$, the only possible source of 
increase of  $\int s^2 d\mu$, $\int |\mathring{r}|^2d\mu$, or $\int |W_+|^2d\mu$
comes from the transition annuli, which already have under control, and so 
we have succeeded in producing a metric $g$ on $\tilde{X}$ with all the claimed properties. 
\end{proof}

Here is a simple application of this Lemma:

\begin{lem}\label{fore}
The $4$-manifold $4 \overline{\bcp}_2$
admits an anorexic sequence. 
\end{lem}
\begin{proof}
Consider the involution of $S^3\times S^1\subset \HH \times \CC$ given by
$(q,z)\mapsto (\bar{q}, \bar{z})$. Now equip $S^3\times S^1$ with the product of the
unit-sphere metric on $S^3$ and the radius-$\epsilon$ metric on $S^1$. These
metrics descend to orbifold metrics on $(S^3\times S^1)/\ZZ_2$ with bounded sectional
curvature, but with arbitrarily small volume; thus we obtain an anorexic sequence
of such metrics by taking any sequence   $\epsilon_j\to 0$. 
The $4$-manifold $\tilde{X}$ obtained from $X$ by replacing its orbifold
singularities with $2$-spheres of self-intersection $-2$ therefore also
admits anorexic sequences by Lemma \ref{scrunch}. But $X=V\# V$,
where $V=S^4/\ZZ_2$, and so $\tilde{X}= \tilde{V}\#\tilde{V}$. However,
as we already noted in the proof of Proposition \ref{change}, $\tilde{V}= 2\overline{\bcp}_2$,
and so $\tilde{X}= 4\overline{\bcp}_2$. 
\end{proof}

We thus obtain our first non-existence result: 

\begin{prop}
The $4$-manifold $4 \overline{\bcp}_2$
does not admit optimal metrics. 
\end{prop}
\begin{proof}
By Lemmata  \ref{skinny} and \ref{fore}, an optimal metric on $4 \overline{\bcp}_2$
would have to be SFASD. However, $(2\chi + 3\tau ) (4 \overline{\bcp}_2)=0$,
so Proposition \ref{fine} would imply that any SFASD metric on the 
simply connected $4$-manifold $4 \overline{\bcp}_2$ would be hyper-K\"ahler.
But such a metric would entail the existence of  non-trivial self-dual harmonic $2$-forms, 
which is excluded here,  since $b_+(4 \overline{\bcp}_2)=0$. 
\end{proof}

A rather  more important application of \ref{scrunch} is the following: 

\begin{lem} \label{niner} The $4$-manifold $\bcp_2\#9 \overline{\bcp}_2$
admits an anorexic sequence. 
\end{lem}
\begin{proof}
Consider the involution  of $S^2\times T^2$ which is  obtained as the product of 
a $180^\circ$ rotation of $S^2$ around an axis and the Weierstrass involution of
an elliptic curve: 
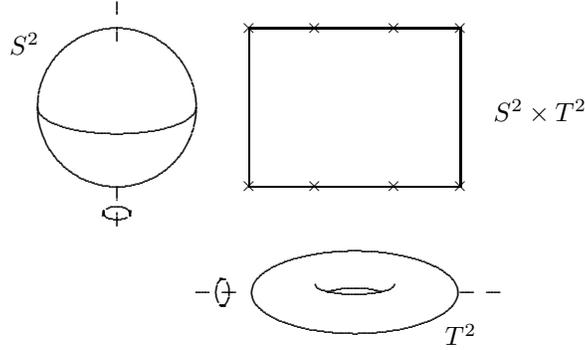
\begin{center}
\mbox{
\beginpicture
\setplotarea x from 0 to 200, y from -5 to 170
\putrectangle corners at 80 130 and 160 70
\put {$_{\times}$} [B1] at 80 70
\put {$_{\times}$} [B1] at 105 70
\put {$_{\times}$} [B1] at 135 70
\put {$_{\times}$} [B1] at 160 70
\put {$_{\times}$} [B1] at 80 130
\put {$_{\times}$} [B1] at 105 130
\put {$_{\times}$} [B1] at 135 130
\put {$_{\times}$} [B1] at 160 130
\put {$T^2$} [B1] at 160 10 
\put {$S^2$} [B1] at -5 120
\circulararc 360 degrees from 30 130 center at 30 100
\ellipticalarc axes ratio 3:1 -180 degrees from 60 100
center at 30 100
\put {$S^2\times T^2$} [B1] at 190 95
\ellipticalarc axes ratio 5:2  360 degrees from 150 40
center at 120 30
\ellipticalarc axes ratio 4:1 -180 degrees from 135 33
center at 120 33
\ellipticalarc axes ratio 4:1 145 degrees from 130 30
center at 120 29
\ellipticalarc axes ratio 1:2 140 degrees from 70 35
center at 70 30
\ellipticalarc axes ratio 1:2 140 degrees from 70 25
center at 70 30
\ellipticalarc axes ratio 2:1 140 degrees from 35 60
center at 30 60
\ellipticalarc axes ratio 2:1 140 degrees from 25 60
center at 30 60
\arrow <2pt> [.1,.3] from 70 35 to 71 35
\arrow <2pt> [.1,.3] from 70 25 to 69 25
\arrow <2pt> [.1,.3] from 35 60 to 35 59
\arrow <2pt> [.1,.3] from 25 60 to 25 61
{\setlinear
\setdashes 
\plot  30 140 30 125  /
\plot  30 70 30 55  /
\plot  60 30 80 30   /
\plot  160 30 180 30   /
}
\endpicture
}
\end{center}
This involution has exactly $8$ fixed points. 
Let $\tilde{X}$ be the manifold which  desingularizes the orbifold $X=[S^2\times T^2]/\ZZ_2$ 
by replacing each of the resulting eight 
singular points with an $S^2$ of self-intersection $-2$. 

Now it is easy to see that $S^2\times T^2$ admits sequences of metrics 
with bounded sectional curvature, but with volume tending to zero: namely, equip 
$T^2$ with a sequence of flat metrics of smaller and smaller area, 
and take the Riemannian product of these metrics with the standard round metric
on $S^2$. Moreover,  the metrics given by this explicit recipe are all
 $\ZZ_2$-invariant, and so give rise to   a sequence of orbifold metrics  
on   $[S^2\times T^2]/\ZZ_2$ with bounded sectional curvature for which the total volume
tends to zero. Such a sequence is  anorexic, and also has  the  special
property that 
$\int |r|^2d\mu \to 0$. By Lemma \ref{scrunch}, $\tilde{X}$ therefore also admits such a special 
anorexic sequence.

It only remains to show  that  $\tilde{X}$ is diffeomorphic to 
$\bcp_2\#9\overline{\bcp}_2$.  To se this,  think of $S^2\times T^2$ as 
 $\bcp_1\times E$, where $E$ is an elliptic curve. Then $\tilde{X}$ becomes a complex surface which has 
 a branched double cover biholomorphic to  
$\bcp_1\times E$ blown up at eight  points. 
This complex surface is simply connected, and it has Kodaira dimension 
$-\infty$ because it contains a   $\bcp_1$ with trivial normal bundle. By the
Enriques-Kodaira classification \cite{bpv}, any such complex surface is rational, and 
hence must be diffeomorphic to either $S^2\times S^2$ or a connected sum
$\bcp_2\# k\overline{\bcp}_2$. However, 
$$(2\chi + 3\tau)(\tilde{X}) = \frac{1}{4\pi^2}\int \left(\frac{s^2}{24}+2|W_+|^2-\frac{|\mathring{r}|^2}{2}\right)d\mu$$
must vanish, since the right-hand side will certainly  tend to zero for
our special anorexic sequence. Since 
\begin{eqnarray*}
(2\chi + 3\tau)(S^2\times S^2) & = & 8 \\
(2\chi + 3\tau)(\bcp_2\# k\overline{\bcp}_2) & = & 9-k , 
\end{eqnarray*}
 it thus follows that 
 $\tilde{X}$ must be 
diffeomorphic to $\bcp_2\# 9\overline{\bcp}_2$. 
\end{proof}

In fact, one does not need to appeal to  any classification machinery to check that 
$\tilde{X}\approx \bcp_2\# 9\overline{\bcp}_2$. Working with one's bare hands \cite{lky},
it is not difficult to show  that, with the fixed complex structure used above, 
$\tilde{X}$ is precisely the complex surface obtained by iteratively blowing up
 $\bcp_2$ at a configuration of points
arranged as in the following diagram, in which
a pair of adjacent points on a line means that one is to blow up the first point, and
then blow up the resulting exceptional divisor at the point corresponding to the direction of
the line, while the adjacent triple of points on a line has an analogous interpretation:
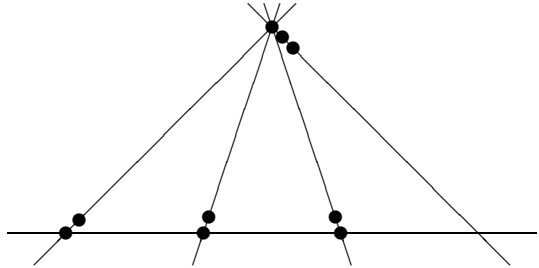
\begin{center}
\begin{picture}(240,120)(0,3)
\put(109,109){\line(-1,-1){99}}
\put(91,109){\line(1,-1){99}}
\put(103,109){\line(-1,-3){33}}
\put(97,109){\line(1,-3){33}}
\put(0,22){\line(1,0){200}}
 \put(100,100){\circle*{5}}
  \put(22,22){\circle*{5}}
  \put(27,27){\circle*{5}}
 \put(104,96){\circle*{5}}
   \put(108,92){\circle*{5}}
 \put(126,22){\circle*{5}}
   \put(124,28){\circle*{5}}
   \put(74,22){\circle*{5}}
   \put(76,28){\circle*{5}} 
\end{picture}
\end{center}
 If we think of the
right-hand line as the line at infinity,   the horizontal line as the $x$-axis, and the 
three other lines as $x=0$, $x=1$ and $x=t$, then the elliptic curves $y^2=Ax(x-1)(x-t)$
foliate an open dense set of the blow-up,  and arise from the $E$ factor 
of $[\bcp_1\times E]/\ZZ_2$. 

\begin{cor}
The simply connected $4$-manifold $\bcp_2 \# 9 \overline{\bcp}_2$ does not admit
optimal metrics. 
\end{cor}
\begin{proof} 
Lemma \ref{niner} allows us to apply the last clause of  Proposition \ref{boney}.  
\end{proof}

To build more complicated examples, first consider the  {\em wormhole space}
obtained by equipping $\RR^4 -\{0\}$ with the metric
$$g_{wh}= \left( 1+\frac{\epsilon}{|\vec{x}|^2}\right)^2 \sum_j (dx^j)^2 .$$
Because $1+{\epsilon}/{|\vec{x}|^2}$ is a harmonic function, this conformally flat 
metric is scalar-flat. However, rewriting this metric in polar coordinates as
$$g_{wh}= \left( \varrho + \frac{\epsilon}{\varrho}\right)^2\left[ 
\frac{d \varrho^2}{\varrho^2} + 
h_{S^3} \right] ,$$
where $h_{S^3}$ is the standard metric on the unit $3$-sphere $S^3$, we immediately
see that there is an isometry $\varrho \mapsto \epsilon /\varrho$ 
of the wormhole which interchanges the two ends of $\RR^4-\{ 0\}\approx \RR \times S^3$.
Thus our wormhole connects two asymptotically Euclidean ends, but has $s\equiv 0$
and $W_+\equiv 0$. Now, on any fixed annulus $\varrho \in [\rho /2, \rho ]$,
the wormhole metric uniformly converges in $C^2$ to the Euclidean metric, and so exactly the same
argument used to glue in Eguchi-Hanson metrics allows us to join two manifolds
with $\int s^2 d\mu < \varepsilon $ and $\int |W_+|^2d\mu < \varepsilon$ by a wormhole neck 
so as to obtain 
a new manifold with $\int s^2 d\mu < 3\varepsilon $ and $\int |W_+|^2d\mu < 3\varepsilon$.
Thus:

\begin{lem}\label{wormhole}
Suppose that $M_1$ and $M_2$ are two smooth compact oriented $4$-manifolds
which admit anorexic sequences. Then their connected sum $M_1\# M_2$
 admits anorexic sequences, too. 
\end{lem}

The final basic building block we will need is the Burns metric \cite{lpa}.
This is an asymptotically flat metric on $\overline{\bcp}_2-\{ p\}$ with
$s\equiv 0$ and $W_+\equiv 0$. Rescaled versions of this  metric, restricted
to the complement of a $\bcp_1$, are explicitly given by 
 $$g_{B, \epsilon}= \frac{d\varrho^2}{1-(\frac{\epsilon}{\varrho})^2}+ \varrho^2\left( \sigma_1^2+\sigma_2^2+ 
\left[1-\left(\frac{\epsilon}{\varrho}\right)^2\right]\sigma_3^2\right)$$
 as a metric on $(\epsilon, \infty ) \times S^3$. As the length-scale parameter
 $\epsilon$ 
tends to zero, we once again get uniform $C^2$ convergence
to the Euclidean metric on any fixed annulus $\varrho \in [\rho / 2, \rho ]$, and
the same gluing argument as before therefore gives us the following result:

\begin{lem}\label{blowup}
Suppose that $M$ is a  smooth compact oriented $4$-manifold
which admits an anorexic sequence. Then $M\# \overline{\bcp}_2$
admits  anorexic sequences, too. 
\end{lem}

We now prove the first part of Theorem \ref{charlie}:

\begin{thm}
Let $j$ and $k$ be  integers such that $j\geq 2$ and $k \geq 9j$. Then 
the simply connected $4$-manifold $j\bcp_2\#k \overline{\bcp}_2$
does not admit optimal metrics. 
\end{thm}
\begin{proof}
By induction on $j$,
Lemmata \ref{niner} and  \ref{wormhole} imply
that the connected sum $j\bcp_2\#9j \overline{\bcp}_2$
of $j$ copies of 
 $\bcp_2\#9 \overline{\bcp}_2$ admits an anorexic sequence. 
Lemma \ref{blowup} and induction on $k-9j$ then gives us the existence of 
 an anorexic sequence on $j\bcp_2\#k \overline{\bcp}_2$ for any $k\geq 9j$, $j\geq 2$. 
 Proposition  \ref{boney} therefore  tells us that there is no optimal metric on any
of these non-spin simply connected $4$-manifolds with $b_+= j \geq 2$. 
\end{proof}

Now let us next 
consider some `exotic' smooth structures on  $\bcp_2\#9\overline{\bcp}_2$.
  Let $X$ denote the 
complex orbifold $[\bcp_1\times E]/\ZZ_2$ used in the proof of Lemma \ref{niner},
and let us observe that we have a  holomorphic orbifold
submersion $X\to \bcp_1/\ZZ_2$ given
by projection to the first factor. Near any non-singular point of $\bcp_1/\ZZ_2$,
this is a locally trivial holomorphic $E$-bundle. Let $p$ and $q$ be two relatively prime integers $\geq 2$, and choose two nonsingular points $a,b\in \bcp_1/\ZZ_2$.
Note that $\bcp_1/\ZZ_2$ is really just a copy of $\bcp_1$ with two marked points
which are considered to be orbifold points of oder $2$. Let us now
also mark the points $a$ and $b$, and consider them to be orbifold points
of order $p$ and $q$. At the same time, we modify $V$ to obtain
a new orbifold $V^0_{p,q}$ by replacing the fiber over $a$ with $E/\ZZ_p$
and the fiber over $b$ with $E/\ZZ_q$, where the two actions are 
generated by translation of $E$ of order $p$ and order $q$, respectively. 
This can be done via a {\em logarithmic transformation} in the sense of 
Kodaira; for example, in a neighborhood of $a$ biholomorphic to the open unit
disk  $D\subset \CC$, we  $E\times D$ with $[E\times D]/\ZZ_p$,
where $\ZZ_p$ acts on $E$ as before, and simultaneously acts on 
$D$ via the action generated by $z\mapsto e^{2\pi i/p}z$. 
We then have a holomorphic  orbifold submersion from $V^0_{p,q}$
to our orbifold $\bcp_1$ with four orbifold points. Now choose the compatible
 flat metric
on $E$ with unit area and  use a partition
of unity to patch the 
a product metric on $V$ with  local product metrics on $[E\times D]/\ZZ_p$
and
$[E\times D]/\ZZ_q$. The result is a Riemannian submersion orbifold metric 
on $V$. If we now scale  down the fiber $E$ with keeping the metric on 
our orbifold $\bcp_1$ fixed, the result is therefore a family of metrics 
on $V^0_{p,q}$ with volume tending to zero while the curvature remains 
uniformly bounded \cite{cheegro,lky}. In particular, this is an anorexic sequence
on $V^0_{p,q}$, and Lemma \ref{scrunch} tells us that the complex surface 
$M^0_{p,q}$
obtain by replacing each of the $8$ orbifold singularities of $V^0_{p,q}$
by $(-2)$-curves also admits anorexic sequences; moreover, 
the manifolds $M^0_{p,q}\# \ell\overline{\bcp}_2$ 
all admit anorexic sequences, too, as a consequence of Lemma \ref{blowup}. 
However,  the manifolds $M_{p,q}$ are the so-called {\em Dolgachev surfaces}.
These Dolgachev surfaces are all homeomorphic to $\bcp_2\# 9\overline{\bcp}_2$,
but Donaldson or Seiberg-Witten invariants  can be used \cite{friedolg} to 
show that no two of them are diffeomorphic. Moreover, the corresponding 
smooth structures on the blow-ups remain distinct, no matter how many times
we blow up, and these smooth structures are moreover all distinct
from the standard one on $\bcp_2 \# k\overline{\bcp}_2$. By 
Proposition \ref{boney}, it follows that none of these smooth manifolds
admits an optimal metric, even though $\bcp_2 \# k\overline{\bcp}_2$
{\em does} admit optimal metrics for $k\geq 14$.

The story is similar for  homotopy $K3$ surfaces. Namely, we can view
$T^4/\ZZ_2$ as an orbifold elliptic fibration over $\bcp_1/\ZZ_2$, and
so modify it by logarithmic transforms of odd order at one fiber. 
The resulting orbifolds $V^1_q$ then admit anorexic sequences as before, 
as do the $4$-manifolds $M^1_q$ obtained by replacing their singular points by 
$2$-spheres of self-intersection $-2$. These manifolds are homeomorphic to 
$K3$ surfaces, but as smooth manifolds they are distinct, not only from  $K3$,
but also 
from each other. Proposition \ref{boney} thus tells us that none
of them admits optimal metrics, even though they are homeomorphic
to $K3$, which {\em does} admit optimal metrics. We have thus proved 
Theorem \ref{baker}: 

\begin{thm} The existence or non-existence of optimal metrics 
depends on the choice of smooth structure. In particular, the 
topological $4$-manifolds $K3$ and $\bcp_2\#k \overline{\bcp}_2$,
$k\geq 14$, admit  infinitely many exotic smooth structures for which 
no  optimal metric exists, even though each also admits a ``standard''
smooth structure for which optimal metrics {\em do} exist. 
\end{thm}
 
 A similar construction \cite{lky} yields  anorexic sequences on many exotic
 manifolds homeomorphic to $j \bcp_2 \# k \overline{\bcp}_2$.  
 Consider the complex orbifold $[\Sigma\times E]/\ZZ_2$, where 
 $\Sigma$ is a hyperelliptic Riemann surface of genus $\mathbf g$:

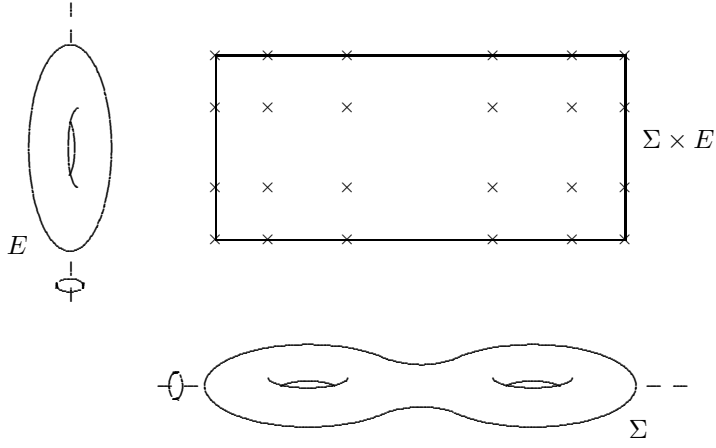
\begin{center}
\mbox{
\beginpicture
\setplotarea x from 0 to 350, y from -5 to 200
\putrectangle corners at 85 155 and 240 85
\put {$\Sigma$} [B1] at 245 10 
\put {$E$} [B1] at  10 80
\put {$\Sigma \times E$} [B1] at  260  120
\put {$_{\times}$} [B1] at 85 105
\put {$_{\times}$} [B1] at 105 105
\put {$_{\times}$} [B1] at 135 105
\put {$_{\times}$} [B1] at 190 105
\put {$_{\times}$} [B1] at 220 105
\put {$_{\times}$} [B1] at 240 105
\put {$_{\times}$} [B1] at 85 85
\put {$_{\times}$} [B1] at 105 85
\put {$_{\times}$} [B1] at 135 85
\put {$_{\times}$} [B1] at 190 85
\put {$_{\times}$} [B1] at 220 85
\put {$_{\times}$} [B1] at 240 85
\put {$_{\times}$} [B1] at 85 135
\put {$_{\times}$} [B1] at 105 135
\put {$_{\times}$} [B1] at 135 135
\put {$_{\times}$} [B1] at 190 135
\put {$_{\times}$} [B1] at 220 135
\put {$_{\times}$} [B1] at 240 135
\put {$_{\times}$} [B1] at 85 155
\put {$_{\times}$} [B1] at 105 155
\put {$_{\times}$} [B1] at 135 155
\put {$_{\times}$} [B1] at 190 155
\put {$_{\times}$} [B1] at 220 155
\put {$_{\times}$} [B1] at 240 155
\ellipticalarc axes ratio 2:5  -360 degrees from 40 150 
center at 30 120 
\ellipticalarc axes ratio 1:4  180 degrees from  33 135
center at  33 120
\ellipticalarc axes ratio 1:4 -145 degrees from  30 130
center at  29 120
\ellipticalarc axes ratio 5:2  280 degrees from 150 40
center at 120 30
\ellipticalarc axes ratio 5:2  -280 degrees from 175 40
center at 205 30
\ellipticalarc axes ratio 4:1 -180 degrees from 135 33
center at 120 33
\ellipticalarc axes ratio 4:1 145 degrees from 130 30
center at 120 29
\ellipticalarc axes ratio 4:1 180 degrees from 190 33
center at 205 33
\ellipticalarc axes ratio 4:1 -145 degrees from 195 30
center at 205 29
\ellipticalarc axes ratio 1:2 140 degrees from 70 35
center at 70 30
\ellipticalarc axes ratio 1:2 140 degrees from 70 25
center at 70 30
\ellipticalarc axes ratio 2:1 140 degrees from 35 68
center at 30 68
\ellipticalarc axes ratio 2:1 140 degrees from 25 68
center at 30 68
\arrow <2pt> [.1,.3] from 70 35 to 71 35
\arrow <2pt> [.1,.3] from 70 25 to 69 25
\arrow <2pt> [.1,.3] from 35 68 to 35 67
\arrow <2pt> [.1,.3] from 25 68 to 25 69
{\setquadratic 
\plot 150 40    163 38    175 40   /
\plot 150 20    163 22    175 20  /
}
{\setlinear
\setdashes 
\plot   27 62  27 82  /
\plot   27 160  27 175  /
\plot  60 30 80 30   /
\plot  245 30 265 30   /
}
\endpicture
}
\end{center}
Giving $E$ its compatible flat metric of various areas,  
products with a fixed metric on $\Sigma$  give us anorexic sequences on the 
orbifold $V^{\mathbf g}= [\Sigma\times E]/\ZZ_2$, and we also get anorexic sequences
on the orbifold $V^{\mathbf g}_q$ obtained by performing a logarithmic
transform to one fiber by the previous gluing argument. If $M^{\bf g}_q$ denotes
the complex surface obtained by replacing the singularities of $V^{\mathbf g}_q$ 
by $(-2)$-curves, then Lemma \ref{scrunch} guarantees that 
$M^{\bf g}_q$ also admits anorexic sequences, and Lemma \ref{blowup}
then tells us that $M^{\bf g}_q\# \ell \overline{\bcp}_2$ admits anorexic sequences,
too. Now $M^{\bf g}_q$ is a simply connected complex surface with
$p_g={\mathbf g}$ and $c_1^2=0$; and it is non-spin if
either  ${\mathbf g}$ or $q$ is even. For any $\ell > 0$, Theorem
\ref{fdmn} tells us that 
$M^{\bf g}_q\# \ell \overline{\bcp}_2$ is  homeomorphic to 
$(2{\mathbf g}+1)\bcp_2 \# (10{\mathbf g}+9+\ell)  \overline{\bcp}_2$;
and we also get the analogous statement for $\ell=0$ if $q$ is even. 
However, by varying $q$, gauge theory can be used to show  \cite{frmo,gost}
 that 
we obtain infinitely many distinct smooth structures in this
way for any fixed $\mathbf g$ and $\ell$.  By Proposition \ref{boney}, 
we therefore have the following result: 

\begin{thm}
For any odd $j\geq 1$ and any $k\geq 5j+4$, the topological manifold
$j\bcp_2 \# k  \overline{\bcp}_2$ admits infinitely distinct smooth 
structures for which no optimal metric exists. 
\end{thm}

Lemma \ref{wormhole} now allows us to produce more examples
of smooth $4$-manifolds with anorexic sequences by 
taking connected sums of our previous examples. Determining whether 
we obtain distinct differentiable structures in this way is a cutting-edge
problem in gauge-theory, however, and it is only by applying the 
sophisticated new machinery of Bauer and Furuta \cite{baufu,bauer2}
that a result of this type can be obtained. Specifically, if ${\mathbf g}$
is {\em odd}, then the bandwidth argument of 
\cite{il1} shows that the 
$M^{\mathbf g}\# M^1_q\# \ell \overline{\bcp}_2$
run through infinitely many differentiable structures on 
$(2{\mathbf g}+4)\bcp_2 \# (10{\mathbf g}
+28+\ell)  \overline{\bcp}_2$ as we vary the even integer  $q\geq 2$, and that 
$2M^{\mathbf g}\# 2M^1_q\# \ell \overline{\bcp}_2$
runs through infinitely many differentiable structures on 
$(4{\mathbf g}+8)\bcp_2 \# (20{\mathbf g}
+56+\ell)  \overline{\bcp}_2$. 
Since Proposition \ref{boney} tells us that none of these spaces 
can admit optimal metrics, we thus obtain the second part of 
 of Theorem 
\ref{charlie}:

\begin{thm}
If $j$ and $k$ are  integers, with $j\geq 5$, $j\not\equiv 0 \bmod 8$,  
and
$k\geq 9j$, then  the topological manifold
$j\bcp_2 \# k  \overline{\bcp}_2$ admits infinitely distinct smooth 
structures for which no optimal metric exists. 
\end{thm}

\section{Concluding Remarks}

While we have seen that many topological $4$-manifolds fail to admit
optimal metrics for many choices of smooth structure, the techniques 
 developed here  do not by any means allow us to determine 
whether or not an optimal metric exists for an {\em arbitrary} smooth structure. 
The reason is that the arguments deployed in \S \ref{nonex}
are
heavily dependent on the existence of anorexic sequences of metrics, 
whereas such sequence simply do not exist for many smooth structures. 
For example, if $X$ is a minimal complex surface of general type, 
and if $M= X\# \ell\overline{\bcp}_2$ is the complex surface
obtained from it by blowing up $\ell$ points, then Seiberg-Witten
theory can be used to show \cite{lric} that any metric on $M$ satisfies
\begin{equation}
\label{west}
\int_M\left(\frac{s^2}{24}+2|W_+|^2\right) d\mu \geq \frac{8\pi^2}{3} c_1^2(X),
\end{equation}
so there are certainly no anorexic metrics on such an $M$. If $X$ is
is simply connected and $\ell$ is sufficiently large, however, this $4$-manifold is 
homeomorphic to 
one of the manifolds treated by Theorem \ref{charlie}, and so represents
an exotic smooth structure on this topological $4$-manifold which is  simply not
amenable to treatment with the current technology. 

The main obstacle to progress on this front  is that 
the estimate (\ref{west}) does not appear to be sharp; if we compare it with known
minimizing sequences for $\int s^2 d\mu$ on a complex surface $M$ of general type 
with minimal
model $X$, we only obtain upper and lower bounds 
$$
\frac{16\pi^2}{3}c_1^2(X) - 8\pi^2 (\chi+ 3\tau )(M) \leq {\mathcal I}_{\mathcal R}(M)
\leq 8\pi^2 c_1^2(X) - 8\pi^2(\chi+ 3\tau )(M)
$$
for ${\mathcal I}_{\mathcal R}(M)$. 
These bounds certainly do not allow one to compute ${\mathcal I}_{\mathcal R}$,
but they do allow us to estimate with sufficient accuracy to be able
to know that at least several different values of ${\mathcal I}_{\mathcal R}$ must
occur for many fixed homeotypes.  An exact computation of ${\mathcal I}_{\mathcal R}$
for such examples  would have many interesting ramifications, and 
should be considered as an outstanding open problem in the subject.

While current technology does not suffice to compute ${\mathcal I}_{\mathcal R}$
for many of the most interesting $4$-manifolds, the  analogous invariants 
\begin{eqnarray*}
	{\mathcal I}_{s}(M) & = & 
 	\inf_{g}\int_{M}|{s}_{g}|^{n/2}d\mu_{g} \\
	{\mathcal I}_{r}(M) & = & 
	\inf_{g}\int_{M}|{r}|^{n/2}_{g}d\mu_{g}  
\end{eqnarray*} 
arising from the scalar and Ricci curvatures {\em do} turn out to be
exactly computable for complex surfaces of general type and many of their
connect sums \cite{lric,il2}.  
It may therefore come as a surprise that 
one key trick  used in computations of  ${\mathcal I}_{r}$
 is closely related to the techniques developed here. 
 Indeed, the Gauss-Bonnet and signature formul{\ae} tell one
that 
$$\int_M |r|^2d\mu = -8\pi^2(2\chi+3\tau )(M) + 8\int_M\left(\frac{s^2}{24}+\frac{1}{2}|W_+|^2\right) d\mu$$
so that the existence of an anorexic sequence is certainly quite sufficient
to allow one to calculate ${\mathcal I}_r$. However, the curious  
difference is that the available Seiberg-Witten 
estimate analogous to (\ref{west})  for this particular combination of $s$ and $W_+$   
 turns out to typically be sharp. For example, if $M=X\#\ell \overline{\bcp}_2$
is a complex surface of general type  with minimal model $X$, one obtains
the estimate
$$\int_M\left(\frac{s^2}{24}+\frac{1}{2}|W_+|^2\right) d\mu \geq 2\pi^2  c_1^2(X),$$
and one can actually find sequences of metrics for which the left-hand side
approaches the expression on the right; thus,  an exact formula 
$$
{\mathcal I}_r(M)= 8\pi^2 [c_1^2(X) + \ell ]
$$
emerges from the discussion. 
Moreover, related argments show that this infimum is unattained whenever 
$\ell >0$. For details  and further applications,   see \cite{lric,il2}.

While Theorem \ref{able} may shed a fair  amount of light on the 
   existence of SFASD metrics on compact $4$-manifolds, 
it by no means  closes the book on the subject. For example, 
we still do not know whether  such metrics exist on $5\overline{\bcp}_2$.
The existence of SFASD metrics on $\bcp_2 \# k \overline{\bcp}_2$ for
$10\leq k\leq 13$ is also not covered  by Theorem \ref{able},
although Yann Rollin and Michael Singer seem to have recently made 
considerable progress concerning these manifolds. 
The global structure of the moduli space of SFASD metrics still remains
a mystery. And it
 would obviously  be of the greatest interest
to sharpen Proposition \ref{class} so as to say something about 
diffeotype when $b_+=0$, or  to 
definitively handle the non-simply connected case. 

We now have many examples of optimal metrics on compact
$4$-manifolds, essentially falling into two classes: the Einstein metrics, and
the scalar-flat anti-self-dual metrics. 
Of course, we can cheaply obtain further examples by reversing 
the orientation of  SFASD manifolds to make them 
self-dual  instead of 
anti-self-dual.  But,   such  trickery aside, 
there do not really seem to be any other known  examples of 
optimal metrics on compact $4$-manifolds. In particular, it  seems
 that all known examples of optimal metrics 
are  critical points of  $\int s^2d\mu$, and so, by optimality, also of 
 $\int |W_+|^2d\mu$. Now a metric on a $4$-manifold is a critical point
 of  $\int |W_+|^2d\mu$ iff it has vanishing  Bach tensor \cite{bes}. Are there any 
optimal metrics on compact $4$-manifolds
that aren't Bach-flat? Are there scalar-flat optimal metrics which 
are Bach-flat, but neither self-dual nor anti-self-dual? Both of these
questions should illustrate the degree to which we still remain
fundamentally ignorant as to  
 the true nature of general optimal metrics, even  
in dimension four.

\bigskip 

{\bf Acknowledgments.}
I would like to express my profound gratitude to  Dominic Joyce, 
who originally convinced me,  by an entirely different line of reasoning,  that
 Theorems \ref{neg} and \ref{pos}  ought to be true. 
 I would also like to warmly thank 
Mike Eastwood and Michael Singer for granting me permission 
 to present some of   their   unpublished results in \S \ref{vanish}.
 Finally, 
I would  like to thank  Chris Bishop, Mark de Cataldo, Sorin Popescu,
and Michael Taylor for helpful suggestions concerning a
number of important   technical details.

  \end{document}